\documentclass[11pt]{article}
\usepackage{amsmath,amsfonts,amssymb,amsthm}
\begin{document}

\newcommand{\1}{{{\bf 1}}}
\newcommand{\id}{{\rm id}}
\newcommand{\Hom}{{\rm Hom}\,}
\newcommand{\End}{{\rm End}\,}
\newcommand{\Res}{{\rm Res}\,}
\newcommand{\Image}{{\rm Im}\,}
\newcommand{\Ind}{{\rm Ind}\,}
\newcommand{\Aut}{{\rm Aut}\,}
\newcommand{\Ker}{{\rm Ker}\,}
\newcommand{\gr}{{\rm gr}}
\newcommand{\Der}{{\rm Der}\,}

\newcommand{\Z}{\mathbb{Z}}
\newcommand{\Q}{\mathbb{Q}}
\newcommand{\C}{\mathbb{C}}
\newcommand{\N}{\mathbb{N}}
\newcommand{\g}{\mathfrak{g}}
\newcommand{\gl}{\mathfrak{gl}}
\newcommand{\h}{\mathfrak{h}}
\newcommand{\wt}{{\rm wt}\;}
\newcommand{\A}{\mathcal{A}}
\newcommand{\D}{\mathcal{D}}
\newcommand{\Lie}{\mathcal{L}}

\def \<{\langle}
\def \>{\rangle}
\def \be{\begin{equation}\label}
\def \ee{\end{equation}}
\def \bex{\begin{exa}\label}
\def \eex{\end{exa}}
\def \bl{\begin{lem}\label}
\def \el{\end{lem}}
\def \bt{\begin{thm}\label}
\def \et{\end{thm}}
\def \bp{\begin{prop}\label}
\def \ep{\end{prop}}
\def \br{\begin{rem}\label}
\def \er{\end{rem}}
\def \bc{\begin{coro}\label}
\def \ec{\end{coro}}
\def \bd{\begin{de}\label}
\def \ed{\end{de}}

\newtheorem{thm}{Theorem}[section]
\newtheorem{prop}[thm]{Proposition}
\newtheorem{coro}[thm]{Corollary}
\newtheorem{conj}[thm]{Conjecture}
\newtheorem{exa}[thm]{Example}
\newtheorem{lem}[thm]{Lemma}
\newtheorem{rem}[thm]{Remark}
\newtheorem{de}[thm]{Definition}
\newtheorem{hy}[thm]{Hypothesis}
\makeatletter \@addtoreset{equation}{section}
\def\theequation{\thesection.\arabic{equation}}
\makeatother \makeatletter

\begin{Large}
\begin{center}
\textbf{On certain vertex algebras and their modules
associated with vertex algebroids}
\end{center}
\end{Large}
\begin{center}{Haisheng Li\footnote{Partially supported by an NSA grant}\\
Department of Mathematical Sciences, Rutgers University, Camden, NJ 08102\\
and\\
Department of Mathematics, Harbin Normal University, Harbin, China\\
\ \
\\
Gaywalee Yamskulna\\
Department of Mathematical Sciences, Binghamton University, NY 13902\\
and\\
Institute of Science, Walailak University, Thailand}
\end{center}

\begin{abstract}
We study the family of vertex algebras associated with vertex
algebroids, constructed by Gorbounov, Malikov, and Schechtman.
As the main result,
we classify all the (graded) simple modules for such vertex
algebras and we show that the equivalence classes of
graded simple modules one-to-one
correspond to the equivalence classes of 
simple modules for the Lie algebroids associated
with the vertex algebroids. To achieve our goal, we construct and
exploit a Lie algebra from a given vertex algebroid.
\end{abstract}

\section{Introduction}
For most of the important examples of vertex operator algebras
$V=\coprod_{n\in \Z}V_{(n)}$ graded by the $L(0)$-weight (see [FLM,
FHL]), the $\Z$-grading satisfies the condition that $V_{(n)}=0$ for
$n<0$ and $V_{(0)}=\C {\bf 1}$ where ${\bf 1}$ is the vacuum vector.
For a vertex operator algebra $V$ with this special property, the
homogeneous subspace $V_{(1)}$ has a natural Lie algebra structure
with $[u,v]=u_{0}v$ for $u,v\in V_{(1)}$ and the product $u_{1}v$
$(\in V_{(0)})$ defines a symmetric invariant bilinear form on
$V_{(1)}$.

In a series of study on
Gerbs of chiral differential operators in [GMS] and on
chiral de Rham complex in [MSV, MS1,2], Malikov and his coauthors investigated
$\N$-graded vertex algebras $V=\coprod_{n\in \N}V_{(n)}$ with $V_{(0)}$
{\em not} necessarily $1$-dimensional.  
In this case, the bilinear operations $(u,v)\mapsto u_{i}v$ for $i\ge 0$
are closed on $V_{(0)}\oplus V_{(1)}$:
$$u_{i}v\in V_{(0)}\oplus V_{(1)} \;\;\;\mbox{ for }u,v\in
V_{(0)}\oplus V_{(1)},\; i\ge 0.$$ The skew symmetry and the Jacobi
identity for the vertex algebra $V$ give rise to several compatibility
relations.  Such algebraic structures on $V_{(0)}\oplus V_{(1)}$ are
summarized in the notion of what was called a $1$-truncated conformal
algebra. Furthermore, the subspace $V_{(0)}$ equipped with the product
$(a,b)\mapsto a_{-1}b$ is a commutative associative algebra with the
vacuum vector ${\bf 1}$ as the identity (cf. [Lia]) and $V_{(0)}$ as a
nonassociative algebra acts on $V_{(1)}$ by $a\cdot u=a_{-1}u$ for
$a\in V_{(0)}$, $u\in V_{(1)}$. All these structures on $V_{(0)}\oplus
V_{(1)}$ are further summarized in the notion of what was called a
vertex $A$-algebroid, where $A$ is a (unital) commutative associative
algebra.  On the other hand, in [GMS], among other important results,
Gorbounov, Malikov and Schechtman constructed an $\N$-graded vertex
algebra $V=\coprod_{n\in \N}V_{(n)}$ {}from any vertex $A$-algebroid,
such that $V_{(0)}=A$ and the vertex $A$-algebroid $V_{(1)}$ is
isomorphic to the given one.  All the constructed $\N$-graded vertex
algebras are generated by $V_{(0)}\oplus V_{(1)}$ with a spanning
property of PBW type.  As it was demonstrated in [GMS], such
$\N$-graded vertex algebras are natural and important to study. For
example, the vertex (operator) algebra associated with a $\beta\gamma$
system, which plays a central role in free field realization of affine
Lie algebras (see [W, FF1-3, FB]) is such an $\N$-graded vertex
algebra.  The vertex (operator) algebras constructed from toroidal Lie
algebras are also of this type (see [BBS, BDT]).

The main purpose of this paper is to classify all the graded simple
modules for the vertex algebras associated with vertex algebroids.
Note that a $1$-truncated conformal algebra is a generalization of a
Lie algebra equipped with a symmetric invariant bilinear form.  For any
Lie algebra $\g$ equipped with a symmetric invariant bilinear form
$\<\cdot,\cdot\>$, one has the affine Lie algebra $\hat{\g}$ (see
[K1]), which is
the one-dimensional (universal) central extension of the loop algebra
$\g\otimes \C[t,t^{-1}]$.  It has been well known ([FZ, Lia],
cf. [Li3, LL]) that for any such given affine Lie algebra $\hat{\g}$
and for any complex number $\ell$, one has a canonical $\N$-graded
vertex algebra $V_{\hat{\g}}(\ell,0)$, where the degree zero subspace
is $\C {\bf 1}$ (one-dimensional) and the degree-one subspace is
naturally identified with $\g$ which generates the whole vertex
algebra $V_{\hat{\g}}(\ell,0)$.  (Here, $V_{\hat{\g}}(\ell,0)$ is a
generalized Verma $\hat{\g}$-module or Weyl module of level $\ell$.)
Furthermore, it was proved (see [Li3, LL]) that the category of
$V_{\hat{\g}}(\ell,0)$-modules is canonically isomorphic to the
category of ``restricted'' $\hat{\g}$-modules of level $\ell$.  
Here, from any $1$-truncated conformal algebra we
construct a Lie algebra, generalizing the construction of affine Lie
algebras.  This Lie algebra is precisely a vertex Lie algebra in the
sense of [DLM]. 
By using this Lie algebra and a result of [DLM], we
construct an $\N$-graded vertex algebra and its modules.  Then
following [GMS], we define the desired vertex algebra $V_{B}$ associated
with a vertex $A$-algebroid $B$ as a quotient vertex algebra.
Associated to each vertex $A$-algebroid $B$ there is a natural
Lie $A$-algebroid (see Section 2 for details). As our main result, 
we prove that the equivalence classes of $\N$-graded simple
$V_{B}$-modules
one-to-one correspond to the equivalence classes of simple modules
for the associated Lie $A$-algebroid.
In [MS2], an analogous result in a lightly different and 
more concrete situation of the chiral de Rham complex had been obtained.
(We are grateful to the referee for pointing out this to us.)

We mention that our construction of the vertex algebra associated with
a $1$-truncated conformal (Lie) algebra is (slightly) different from
the one given in [GMS].  In [GMS], Gorbounov, Malikov and
Schechtman first constructed a conformal algebra in the sense of
[K2] (where the details of the proof were left as an exercise) and
then they gave a new construction of vertex algebras from conformal
algebras (cf. [K2, P]).  On the other hand, in our construction,
the main work is to construct the Lie algebra, which is a vertex Lie
algebra in the sense of [DLM], from a $1$-truncated conformal algebra.
It is well known that
the affinization of a conformal (Lie) algebra gives rise to
a Lie algebra (see [K2, P]), just as 
the affinization of a vertex algebra gives rise to
a Lie algebra (see [B, FFR, Li1,2, MP]).
The Lie algebra we constructed is precisely the Lie algebra of
the conformal (Lie) algebra associated to the $1$-truncated conformal
algebra by [GMS]. However, 
it seems that completing the exercise posted in [GMS] to fulfill
the details in the construction of the conformal algebra from a $1$-truncated
conformal algebra is
not easier than constructing the Lie algebra. In fact,
we only foresee a complete solution to the exercise by
using our Lie algebra result together with a result of Primc [P].
Notice that in constructing a vertex algebra from a conformal
algebra $C$, Gorbounov, Malikov and Schechtman also exploited a Lie
algebra structure on $C$, which is a proper Lie subalgebra of what we
exploited.

This paper is organized as follows: In Section 2 we recall 
the notions of vertex $A$-algebroid and $1$-truncated
conformal algebra. In Section 3, we construct a Lie algebra from any
$1$-truncated conformal algebra and then we construct a vertex
algebra.  In Section 4, for a given vertex $A$-algebroid we construct
an $\N$-graded vertex algebra and we classify all the
graded simple modules.

\section{$1$-truncated conformal algebras, vertex algebroids and
 Lie algebroids}
In this section, we begin with reviewing some basic notations,
formulas and properties for vertex algebras, 
which will be used in later sections.
We then recall the definitions of the notions of 
$1$-truncated conformal algebra, vertex algebroid and Lie
algebroid.
We also present certain basic properties and
relations among these algebraic objects.

First, we recall the definition of a vertex algebra (cf. [LL]).
A {\em vertex algebra} is a vector space $V$ equipped with a linear map 
\begin{eqnarray}
Y: V&\rightarrow& (\End V)[[x,x^{-1}]],\nonumber\\ 
v&\mapsto& Y(v,x)=\sum_{n\in\Z}v_{n}x^{-n-1}
\;\;\;(\mbox{where }v_{n}\in \End V)
\end{eqnarray}
and equipped with a distinguished vector ${\bf 1}$, 
called the {\em vacuum (vector)},
such that for $u,v\in V$,
\begin{eqnarray}
& &u_{n}v=0\;\;\;\mbox{ for n sufficiently large},\\
& &Y({\bf 1},x)=1,\\
& &Y(v,x){\bf 1}\in V[[x]]\;\;\mbox{ and } 
\lim_{x\rightarrow 0}Y(v,x){\bf 1}=v
\end{eqnarray}
and such that
\begin{eqnarray}
& &x_{0}^{-1}\delta\left(\frac{x_{1}-x_{2}}{x_{0}}\right)Y(u,x_{1})Y(v,x_{2})
-x_{0}^{-1}\delta\left(\frac{x_{2}-x_{1}}{-x_{0}}\right)Y(v,x_{2})Y(u,x_{1})\nonumber\\
& &\ \ \ \ =x_{2}^{-1}\delta\left(\frac{x_{1}-x_{0}}{x_{2}}\right)Y(Y(u,x_{0})v,x_{2})
\end{eqnarray}
the {\em Jacobi identity}.

A {\em $V$-module} is a vector space $W$ equipped with a linear map $Y_{W}$
{} from $V$ to $(\End W)[[x,x^{-1}]]$, where $Y_{W}(v,x)=\sum_{n\in \Z}v_{n}x^{-n-1}$ for $v\in V$,
such that for $v\in V,\; w\in W$,
\begin{eqnarray}
& &v_{n}w=0\;\;\;\mbox{ for n sufficiently large},\\
& &Y_{W}({\bf 1},x)=1,\\
& &x_{0}^{-1}\delta\left(\frac{x_{1}-x_{2}}{x_{0}}\right)Y_{W}(u,x_{1})Y_{W}(v,x_{2})
-x_{0}^{-1}\delta\left(\frac{x_{2}-x_{1}}{-x_{0}}\right)Y_{W}(v,x_{2})Y_{W}(u,x_{1})\nonumber\\
& &\ \ \ \ =x_{2}^{-1}\delta\left(\frac{x_{1}-x_{0}}{x_{2}}\right)Y_{W}(Y(u,x_{0})v,x_{2})
\end{eqnarray}
for $u,v\in V$.

{}From the Jacobi identity we have Borcherds'
commutator formula and iterate formula:
\begin{eqnarray}
[u_{m},v_{n}]&=&\sum_{i\ge 0}\binom{m}{i}(u_{i}v)_{m+n-i},\\
(u_{m}v)_{n}w&=&\sum_{i\ge 0}(-1)^{i}\binom{m}{i}
\left( u_{m-i}v_{n+i}w-(-1)^{m}v_{m+n-i}u_{i}w\right)
\end{eqnarray}
for $u,v,w\in V,\; m,n\in \Z$.
Define a linear operator $\mathcal{D}$ on $V$ by
\begin{equation*}
{\mathcal{D}}(v)=v_{-2}\1\;\;\;\mbox{ for }v\in V.
\end{equation*}
Then
\begin{eqnarray}
Y(v,x){\bf 1}=e^{x\D}v\;\;\;\mbox{ for }v\in V.
\end{eqnarray}
Furthermore,
\begin{eqnarray}\label{ed-bracket-formula}
[\D,v_{n}]=(\D v)_{n}=-n v_{n-1}
\end{eqnarray}
for $v\in V,\; n\in \Z$. 

%Let $S$ be a subset of a vertex algebra $V$. We say
%that {\em $S$ generates $V$} if 
%$$V={\rm span}\{ u^(1)_{n_1}...u^(r)_{n_r}{\bf 1}\ \ |\ \  r\in
%\N,\ \ u^1,...,u^r\in S,\ \  n_1,...,n_r\in \Z\}.$$

%An {\em ideal} of a vertex algebra $V$ is a subspace $I$
%such that for all $v\in V$, $w\in I$, $\ \ Y(v,x)w\in
%I[[x,x^{-1}]],\text{ and }\ \  Y(w,x)v\in I[[x,x^{-1}]].$

%A subspace $I$ of the vertex algebra $V$ is an ideal if $\D
%I\subset I$, and $v_nw\in I$ for all $v\in V, w\in I, n\in \Z$.
%For a subset $S$ of the vertex algebra $V$,
%$$\<S\>={\rm span}\{v_n\D^i(u)\ \ |\ \ v_n\in V,\ \ n\in \Z,\ \ i\geq
%0,\ \ u\in S\}.$$

A vertex algebra $V$ equipped with a $\Z$-grading $V=\coprod_{n\in
\Z}V_{(n)}$ is called a {\em $\Z$-graded vertex algebra} if ${\bf
1}\in V_{(0)}$ and if for $u\in V_{(k)}$ with $k\in \Z$ and for
$m,n\in \Z$,
\begin{eqnarray}
u_{m}V_{(n)}\subset V_{(k+n-m-1)}.
\end{eqnarray}
An $\N$-graded vertex algebra is defined in the obvious way.

The following notion of $1$-truncated conformal
algebra is due to [GMS]:

\bd{dtruncatedCA} {\em A {\em $1$-truncated conformal
algebra} is a graded vector space $C=C_{0}\oplus C_{1}$, equipped
with a linear map $\partial: C_{0}\rightarrow C_{1}$ and bilinear
operations $(u,v)\mapsto u_{i}v$ for $i=0,1$ of degree $-i-1$ 
on $C=C_{0}\oplus C_{1}$ such that the
following axioms hold:

(Derivation) for $a\in C_{0},\; u\in C_{1}$,
\begin{eqnarray}\label{e2.7}
(\partial a)_{0}=0;\;\; (\partial a)_{1}=-a_{0};\;\;
\partial (u_{0}a)=u_{0}\partial a
\end{eqnarray}

(Commutativity) for $a\in C_{0},\; u,v\in C_{1}$,
\begin{eqnarray}\label{e2.8}
u_{0}a=-a_{0}u;\;\; u_{0}v=-v_{0}u+\partial (v_{1}u);\;\;
u_{1}v=v_{1}u
\end{eqnarray}

(Associativity) for $\alpha,\beta,\gamma\in C_{0}\oplus C_{1}$,
\begin{eqnarray}\label{e2.9}
\alpha_{0}\beta_{i}\gamma
=\beta_{i}\alpha_{0}\gamma +(\alpha_{0}\beta)_{i}\gamma.
\end{eqnarray}}
\ed

\bex{rtruncatedCA} {\em Let $V=\coprod_{n\in \N}V_{(n)}$ be an
$\N$-graded vertex algebra. One can easily show that
$V_{(0)}\oplus V_{(1)}$ is naturally a $1$-truncated conformal algebra. If
$C=\coprod_{n\in \N}C_{(n)}$ is an $\N$-graded conformal algebra
in the sense of [K2], then it is straightforward to see that
$C_{(0)}\oplus C_{(1)}$ is a $1$-truncated conformal algebra. }
\eex

\bex{eample0} 
{\em Let $\g$ be a Lie algebra over $\C$ equipped with a
symmetric invariant bilinear form $\<\cdot,\cdot\>$.
Let $\partial$ be the zero map from $\C$ to $\g$. Then
$C=\C\oplus \g$ is a 1-truncated conformal algebra.}
 \eex

\br{exam2} 
{\em Just as with Lie algebras, direct sums of 
$1$-truncated conformal algebras 
are $1$-truncated conformal algebras.}
\er

\bl{lsimplify}
In the definition of the notion of $1$-truncated conformal algebra,
the third relation in the derivation axiom follows from
the other derivation relations and commutativity relations.
In the associativity axiom, it is sufficient just to assume
those for $\alpha,\beta,\gamma\in C_{1}$ with $i=0,1$ and
those for $\alpha\in C_{0},\;\beta,\gamma\in C_{1}$ with $i=0$.
\el

\begin{proof} For $a\in C_{0},\; u\in C_{1}$, we have
\begin{eqnarray*}
u_{0}\partial a=-(\partial a)_{0}u+\partial ((\partial a)_{1}u)
=\partial ((\partial a)_{1}u)=-\partial (a_{0}u)=\partial (u_{0} a),
\end{eqnarray*}
proving the first assertion.  By linearity, we may assume
$\alpha,\beta,\gamma\in C_{0}\cup C_{1}$ in the associativity
axiom. If more than one of $\alpha,\beta,\gamma$ are elements of
$C_{0}$, or if one of $\alpha,\beta,\gamma$ are elements of $C_{0}$ with
$i=1$, due to the degree assumption we have
$$\alpha_{0}\beta_{i}\gamma=\beta_{i}\alpha_{0}\gamma
=(\alpha_{0}\beta)_{i}\gamma=0.$$
Then the nontrivial associativity relations are
those for $\alpha,\beta,\gamma\in C_{1}$ 
and those for $\alpha,\beta,\gamma\in C_{0}\cup C_{1}$ 
with one of them from $C_{0}$ and with $i=0$.
Assume that the associativity relation with 
$\alpha\in C_{0},\;\beta,\gamma\in C_{1}$ with $i=0$ holds.
Using the commutativity relation $a_{0}u=-u_{0}a$ 
for $a\in C_{0},\; u\in C_{1}$, we get
\begin{eqnarray*}
& &\beta_{0}\alpha_{0}\gamma=\alpha_{0}\beta_{0}\gamma-(\alpha_{0}\beta)_{0}\gamma=
\alpha_{0}\beta_{0}+(\beta_{0}\alpha)_{0}\gamma,\\
& &\beta_{0}\gamma_{0}\alpha=-\beta_{0}\alpha_{0}\gamma
=-\alpha_{0}\beta_{0}\gamma-(\beta_{0}\alpha_{0})_{0}\gamma
=(\beta_{0}\gamma)_{0}\alpha+\gamma_{0}\beta_{0}\alpha.
\end{eqnarray*}
This proves the second assertion.
\end{proof}

{\em A Leibniz algebra} is a nonassociative algebra
$\Gamma$ satisfying the following condition:
\begin{eqnarray}
u\cdot (v\cdot w)=(u\cdot v)\cdot w+v\cdot (u\cdot w)\;\;\;\mbox{
for } u,v,w\in \Gamma.
\end{eqnarray}
Note that any Lie algebra is a Leibniz algebra. In particular,
for any vector space $W$ the general linear Lie algebra $\gl(W)$ is
a Leibniz algebra.
A {\em representation of Leibniz algebra} $\Gamma$ 
on a vector space $W$ is
a Leibniz algebra homomorphism $\rho$ from $\Gamma$ to $gl(W)$, i.e.,
a linear map from $\Gamma$ to $\End W$ such that
\begin{eqnarray}
\rho(u\cdot v)=[\rho(u),\rho(v)]\;\;\;\mbox{ for } u,v\in \Gamma.
\end{eqnarray}
Just as with a Lie algebra, for any two $\Gamma$-modules $U$ and
$V$, we have a $\Gamma$-module $U\oplus V$.

In terms of these notions we have (cf. [GMS]):

\bp{psummary} 
Let $C=C_0\oplus C_1$ be a graded vector space
equipped with a linear map $\partial$ from $C_{0}$ to $C_{1}$
and equipped with bilinear maps $(u,v)\mapsto u_{i}v$ of degree $-i-1$
on $C=C_{0}\oplus C_{1}$ for $i=0,1$.
Then $C$ is a $1$-truncated conformal algebra if and only if
\begin{enumerate}
\item The pair $(C_1,[\cdot,\cdot])$ carries the structure of
a Leibniz algebra where $[u,v]=u_0v$ for $u,v\in C_1$. 
\item The space $C_0$ is a $C_1$-module with
$u\cdot a=u_0a$ for $u\in C_1, a\in C_0$. 
\item The map $\partial$ is a $C_1$-module homomorphism. 
\item The subspace $\partial C_0$ of $C_{1}$ annihilates the $C_1$-module 
$C_0\oplus C_1$. 
\item The bilinear linear map $\<\cdot,\cdot\>$ from 
$C_1\otimes C_1$ to $C_0$
defined by $\<u,v\>=u_1v$ for $u,v\in C_1$ is 
a $C_{1}$-module homomorphism and furthermore
\begin{eqnarray}
& &u_{0}a=-a_{0}u,\label{e2.11}\\
& &\<\partial a,u\>=-a_{0}u,\label{e2.12}\\
& &[u,v]+[v,u]=\partial\<u,v\>\label{e2.13}\\
& &\<u,v\>=\<v,u\>\label{e2.14}
\end{eqnarray}
for $a\in C_{0},\; u,v\in C_{0}$.
\end{enumerate}
\ep

\begin{proof} 
Notice that the derivation axiom (\ref{e2.7}) amounts to
Property 4, (\ref{e2.12}) and Property 3, and that
the commutativity axiom (\ref{e2.8}) amounts to
(\ref{e2.11}), (\ref{e2.13}) and (\ref{e2.14}).

The associativity axiom (\ref{e2.9}) 
with $\alpha,\beta,\gamma\in C_{1},\; i=0$
amounts to Property 1 stating that $C_1$ is a Leibniz algebra
and  the associativity axiom (\ref{e2.9}) 
with $\alpha,\beta\in C_{1},\;\gamma\in C_{0},\; i=0$
amounts to Property 2 stating that $C_0$ is a $C_1$-module.
{}From the proof of Lemma \ref{lsimplify}, we see that
the associativity relation (\ref{e2.9}) 
with $\alpha,\beta\in C_{1},\;\gamma\in C_{0},\; i=0$
implies the associativity relation (\ref{e2.9}) 
with $\alpha\in C_{0},\;\beta,\gamma\in C_{1},\; i=0$.
Furthermore, the associativity axiom (\ref{e2.9}) 
with $\alpha,\beta,\gamma\in C_{1},\; i=1$
amounts to Property 5. Now our assertion follows from Lemma \ref{lsimplify}.
\end{proof}

As an immediate consequence we have:

\bc{cLie}
Let $C=C_{0}\oplus C_{1}$ be a $1$-truncated conformal algebra.
Then the quotient algebra $C_1/\partial C_0$ of the Leibniz algebra $C_{1}$ 
is a Lie algebra. Moreover, $C_0$ is a $C_{1}/C_{0}$-module
with $(u+\partial a)\cdot b=u_{0}b$ for $u\in C_{1},\; a,b\in C_{0}$. 
\ec

Next, we recall the following definition of a vertex algebroid from
[GMS, Br1-2]:

\bd{dalgebroid}
{\em Let $A$ be a (unital) commutative associative algebra (over $\C$).
A {\em vertex $A$-algebroid} is a $\C$-vector space $\Gamma$ equipped with

(0) a $\C$-bilinear map
\begin{eqnarray}
A\times \Gamma \rightarrow \Gamma;\;\;\;
(a,v)\mapsto a*v
\end{eqnarray}
such that $1*v=v$ (i.e., a ``non-associative unital $A$-module'')

(1) a structure of a Leibniz $\C$-algebra
$[\cdot,\cdot]:
    \Gamma\otimes_{\C}\Gamma\rightarrow \Gamma$

(2) a homomorphism of Leibniz $\C$-algebras
     $\pi: \Gamma\rightarrow {\rm Der}(A)$

(3) a symmetric $\C$-bilinear pairing
$\<\cdot,\cdot\>: \Gamma\otimes_{\C}\Gamma\rightarrow A$

(4) a $\C$-linear map $\partial: A\rightarrow \Gamma$ such that
    $\pi\circ \partial =0$,\\
which satisfy the following conditions:
\begin{eqnarray}
a*(a'*v)-(aa')*v&=&\pi(v)(a)* \partial (a')+\pi(v)(a')* \partial (a)
\label{pe2.22}\\
{[u,a*v]}&=&\pi(u)(a)*v+a*[u,v]\label{pe2.23}\\
{[u,v]+[v,u]}&=& \partial(\<u,v\>)\label{e2.23}\\
\pi (a*v)&=& a\pi(v)\label{epi-hom}\\
\<a*u,v\>&=& a\<u,v\>-\pi(u)(\pi(v)(a))\label{e2.26}\\
\pi(v)(\<v_{1},v_{2}\>)&=&\<[v,v_{1}],v_{2}\>+\<v_{1},[v,v_{2}]\>
\label{esymmetricmap}\\
\partial (aa')&=& a* \partial (a')+a'* \partial (a)\label{epartial-der}\\
{[v,\partial (a)]}&=& \partial (\pi(v)(a))\label{epartial-hom}\\
\<v,\partial (a)\>&=& \pi(v)(a)\label{e2.29}
\end{eqnarray}
for $a,a'\in A,\; u, v,v_{1}, v_{2}\in \Gamma$.}
\ed

\br{rvertex-alg}
{\em  Let $V=\coprod_{n\in \N}V_{(n)}$ be an $\N$-graded vertex algebra. 
Then $A=V_{(0)}$ is the commutative associative
algebra with ${\bf 1}$ as identity where $aa'=a_{-1}a'$ for $a,a'\in A$ 
and $V_{(1)}$ is a vertex $A$-algebroid (see [MS1]).}
\er

\br{rtypo}
{\em In the original definition given in [Br1-2] 
the two terms on the left hand side of (\ref{pe2.22})
have a negative sign. According to the calculation for Remark
\ref{rvertex-alg} from an $\N$-graded vertex algebra,
the current (\ref{pe2.22}) is the correct one.}
\er

\begin{prop}\label{pact} 
Let $A$ be a unital commutative associative algebra
and $B$ a module for $A$ as a nonassociative algebra.
Then a vertex $A$-algebroid structure on $B$ exactly amounts to
a $1$-truncated conformal algebra structure on $C=A\oplus B$ with
\begin{eqnarray}
& &a_{i}a'=0,\\
& &u_{0}v=[u,v],\;\; \ \ u_{1}v=\<u,v\>,\\
& &u_{0}a=\pi(u)(a),\;\; \ \ a_{0}u=-u_{0}a=-\pi(u)(a)
\end{eqnarray}
for $a,a'\in A,\; u,v\in B,\; i=0,1$, such that
\begin{eqnarray}
& &a(a'u)-(aa')u=(u_{0}a)\partial a'+(u_{0}a')\partial a,\label{eabu-abu}\\
& &u_{0}(av)-a(u_{0}v)=(u_{0}a)v,\label{eu0av2.36}\\
& &u_{0}(aa')=a(u_{0}a')+(u_{0}a)a',\label{u0ab2.35}\\
& &a_{0}(a'v)=a'(a_{0}v),\label{ea0a'v}\\
& &(au)_{1}v=a(u_{1}v)-u_{0}v_{0}a,\label{eau1v}\\
& &\partial (aa')=a\partial a'+a'\partial a.\label{epab}
\end{eqnarray}
\end{prop}

\begin{proof} 
We shall just compare the properties (\ref{epab})--(\ref{eabu-abu}) 
and those for 1-truncated conformal algebras in Proposition
\ref{psummary} with the axioms in Definition \ref{dalgebroid}.

First, Property 1 in Proposition \ref{psummary} corresponds to Axiom
(1) in Definition \ref{dalgebroid}, and the relation (\ref{e2.14}) corresponds to Axiom (3) in Definition \ref{dalgebroid}.
Second, the relation (\ref{u0ab2.35}) and Property 2 in Proposition \ref{psummary}
correspond to Axiom (2) in Definition \ref{dalgebroid}
stating that $\pi: \Gamma \rightarrow  \Der (A)$ is
a homomorphism of Leibniz algebras.
Third,  Property 3, the first part of Property 5 and
(\ref{e2.13}) in Proposition \ref{psummary}
correspond to (\ref{epartial-hom}), (\ref{esymmetricmap}) and
(\ref{e2.23}), respectively.

Note that the relation (\ref{e2.11}) is used in defining $a_{0}u$.
In the presence of this relation and (\ref{e2.14}) (namely, Axiom (3)),
(\ref{e2.29}) amounts to (\ref{e2.12}).
Due to the commutativity relation $a_{0}v=-v_{0}a$
for $a\in A,\; v\in B$, we see that (\ref{ea0a'v})
is equivalent to that
$(av)_{0}a'=a(v_{0}a')$ for $a,a'\in A,\; v\in B$, which corresponds
to (\ref{epi-hom}).

Finally, Property 4 in Proposition \ref{psummary} implies Axiom (4) in Definition \ref{dalgebroid}.
On the other hand, Axioms (3), (4), conditions (\ref{e2.23}), (\ref{epartial-hom}), and (\ref{e2.29})
in Definition \ref{dalgebroid} imply Property 4 in Proposition \ref{psummary}. Clearly, (\ref{eabu-abu}), (\ref{epab}), (\ref{e2.26}) and (\ref{eu0av2.36})
correspond to (\ref{pe2.22}), (\ref{epartial-der}), (\ref{eau1v}) and
(\ref{pe2.23}), respectively.
\end{proof}

\begin{de}
{\em Let $A$ be a commutative associative algebra.
A {\em Lie $A$-algebroid} is a Lie algebra $\g$ equipped
with an $A$-module structure and a module action on $A$ by
derivation such that
\begin{eqnarray}
[u,av]&=&a[u,v]+(ua)v,\\
a(ub)&=&(au)b\;\;\;\mbox{ for }u,v\in \g,\; a,b\in A.
\end{eqnarray}

A {\em module} for a Lie $A$-algebroid $\g$ is a vector space $W$
equipped with a $\g$-module structure and an $A$-module structure
such that
\begin{eqnarray}
& &u(aw)-a(uw)=(ua)w,\\
& &a(uw)=(au)w\;\;\;\mbox{ for }a\in A,\; u\in \g,\; w\in W.
\end{eqnarray}}
\end{de}

The following is an immediate consequence (see [Br2]):

\bl{verlie} 
Let $A$ be a commutative associative algebra (over
$\C$) and let $B$ be a vertex $A$-algebroid. Then $B/A\partial A$ is
a Lie $A$-algebroid.
\el

\br{rlie-vertex-broid} 
{\em Let $A$ be a unital commutative associative algebra and
let $B$ be a Lie $A$-algebroid such that $B$ acts as zero on $A$.
Let $\<\cdot, \cdot\>$ be the zero map from $B\otimes_{\C} B$ to $A$
and let $\partial$ be the zero map from $A$ to $B$. Then $B$ is a
vertex $A$-algebroid.}
\er

Let $\g$ be a Lie $A$-algebroid. Then $\g$ is a Lie algebra
with $A$ a $\g$-module. By adjoining the $\g$-module $A$ to $\g$ we have a
Lie algebra $A\oplus \g$ with $A$ as an abelian ideal.
Denote by $J$ the 2-sided ideal of the universal enveloping algebra
$U(A\oplus \g)$ generated by the vectors
\begin{eqnarray}
e-1,\ \  a\cdot a'-aa', \ \ a\cdot u-au
\end{eqnarray}
for $a, a'\in A,\; u\in \g$, where $\cdot$ denotes the product
in the universal enveloping algebra. Set
\begin{eqnarray}
\bar{U}(A\oplus \g)=U(A\oplus \g)/J.
\end{eqnarray}

\bp{pkey} Let $\g$ be a Lie $A$-algebroid.
Then a (simple) module structure for the Lie
$A$-algebroid $\g$ 
on a vector space $W$ exactly amounts to 
a (simple) $\bar{U}(A\oplus \g)$-module 
structure.
\ep

\begin{proof} From definition, a module structure 
for the Lie $A$-algebroid $\g$ on $W$ amounts to
a linear map $\psi$ from $A\oplus \g$ to $\End W$ such that
\begin{eqnarray*}
& &\psi(e)=1,\ \ \ 
\psi(aa')=\psi(a)\psi(a'),\\
& &\psi([u,v])=\psi(u)\psi(v)-\psi(v)\psi(u),\\
& &\psi(au)=\psi(a)\psi(u),\\
& &\psi(u)\psi(a)-\psi(a)\psi(u)=\psi(u(a))
\end{eqnarray*}
for $a,a'\in A,\; u,v\in \g$. On the other hand, such a linear map $\psi$
exactly amounts to a $\bar{U}(A\oplus \g)$-module structure on $W$.
Then it is clear.
\end{proof}

\section{Lie algebras and vertex algebras associated 
with $1$-truncated conformal algebras}
In this section we construct an honest Lie algebra from any
$1$-truncated conformal algebra.  Then using
this Lie algebra and a result of [DLM] we construct a vertex algebra for any
$1$-truncated conformal algebra.
The constructed Lie algebra is precisely a vertex Lie algebra in the
sense of [DLM] and the construction generalizes the
loop construction of affine Lie algebras from a Lie algebra equipped
with a symmetric invariant bilinear form. 

We start with a graded vector space $C=A\oplus B$ equipped with the following
linear maps:
\begin{eqnarray}
\partial: A&\rightarrow& B;\; a\mapsto \partial a,\label{emap-D}\\
A\times B&\rightarrow& A; \;(a, b)\mapsto a_{0}b,\\
B\times A&\rightarrow& A; \;(b,a)\mapsto b_{0}a,\\
B\times B&\rightarrow& B; \;(b, b')\mapsto b_{0}b',\\
B\times B&\rightarrow& A; \;(b, b')\mapsto
b_{1}b'.\label{emap-b1b}
\end{eqnarray}
Set
\begin{eqnarray}
L(A\oplus B)=(A\oplus B)\otimes \C[t,t^{-1}],
\end{eqnarray}
a vector space with subspaces $L(A)$ and $L(B)$, 
which are defined in the obvious way.
Furthermore, set
\begin{eqnarray}
\hat{\partial}=\partial\otimes 1+1\otimes d/dt: 
\;\; L(A)\rightarrow L(A\oplus B).
\end{eqnarray}
Define
\begin{eqnarray}
\deg (a\otimes t^{n})&=&-n-1\;\;\;\mbox{ for }a\in A,\;n\in \Z,\\
\deg (b\otimes t^{n})&=&-n\;\;\;\mbox{ for }b\in B,\;n\in \Z.
\end{eqnarray}
Then $L(A\oplus B)$ becomes a $\Z$-graded vector space:
\begin{eqnarray}\label{eZ-grading}
L(A\oplus B)=\coprod_{n\in \Z}L(A\oplus B)_{(n)},
\end{eqnarray}
where
\begin{eqnarray*}
L(A\oplus B)_{(n)}=A\otimes \C t^{-n-1}+ B\otimes \C t^{-n}.
\end{eqnarray*}
The subspaces $L(A)$ and $L(B)$ are graded subspaces and for $n\in \Z$ we have
\begin{eqnarray*}
L(A)_{(n)}=A\otimes \C t^{-n-1}=\{ a\otimes t^{-n-1}\;|\; a\in A\}.
\end{eqnarray*}
The linear map $\hat{\partial}$ (from $L(A)$ to $L(A\oplus B)$) is homogeneous of degree $1$
and we have
\begin{eqnarray}\label{ehatDla}
(\hat{\partial}L(A))_{(m)}=\hat{\partial}L(A)_{(m-1)}
=\{ \hat{\partial}(a\otimes t^{-m})\;|\; a\in A\}
\end{eqnarray}
for $m\in \Z$.

Define a bilinear product ``$[\cdot,\cdot]$'' on $L(A\oplus B)$ by
\begin{eqnarray}
[a\otimes t^{m}, a'\otimes t^{n}]&=&0,\label{edef-bracket-1}\\
{[a\otimes t^{m}, b\otimes t^{n}]}&=&a_{0}b\otimes t^{m+n},\label{edef-bracket-2}\\
{[b\otimes t^{n},a\otimes t^{m}]}&=&b_{0}a\otimes
t^{m+n},\label{edef-bracket-22}\\
{[b\otimes t^{m}, b'\otimes t^{n}]}&=&b_{0}b'\otimes
t^{m+n}+m(b_{1}b')\otimes t^{m+n-1}\label{edef-bracket-3}
\end{eqnarray}
for $a,a'\in A,\; b,b'\in B,\;m,n\in \Z$.
Clearly, the product $[\cdot,\cdot]$ is homogeneous of degree zero.
Now, we have a $\Z$-graded nonassociative algebra $L(A\oplus B)$,
equipped with the product $[\cdot,\cdot]$.

Next, we determine the conditions under which
the product $[\cdot,\cdot]$ reduces to a Lie algebra structure 
on the quotient space $L(A\oplus B)/\hat{\partial}L(A)$.

First we have:

\bp{pto-quotient}
The subspace $\hat{\partial}L(A)$ of the nonassociative algebra
$(L(A\oplus B),[\cdot,\cdot])$ is a two-sided ideal if and only if
\begin{eqnarray}
a_{0}\partial a'&=&0\;\;\;\mbox{ for }a,a'\in A,\label{equotient-1}\\
b_{0}\partial a'&=&\partial b_{0}a'\;\;\;\mbox{ for }b\in B,\; a'\in A,\label{equotient-2}\\
b_{1}\partial a'&=&b_{0}a'\;\;\;\mbox{ for }b\in B,\; a'\in A,\label{equotient-2more}\\
(\partial a)_{0}a'&=&0\;\;\;\mbox{ for }a,a'\in A,\label{equotient-1new}\\
(\partial a)_{0}b&=&0\;\;\;\mbox{ for }a\in A,\;b\in B,\label{equotient-3}\\
(\partial a)_{1}b&=&-a_{0}b\;\;\;\mbox{ for }a\in A,\;b\in
B.\label{equotient-4}
\end{eqnarray}
\ep

\begin{proof} For $a,a'\in A,\; m,n\in \Z$,
by (\ref{edef-bracket-1})--(\ref{edef-bracket-22}) we have
\begin{eqnarray}\label{eaDa}
[a\otimes t^{m}, \hat{\partial}(a'\otimes t^{n})]&=&[a\otimes t^{m},
\partial a'\otimes t^{n}+na'\otimes t^{n-1}]\nonumber\\
&=&a_{0}\partial a'\otimes t^{m+n},
\end{eqnarray}
\begin{eqnarray}\label{eDa-a}
[\hat{\partial}(a'\otimes t^{n}),a\otimes t^{m}]&=&[ \partial a'\otimes
t^{n}+na'\otimes t^{n-1},a\otimes
t^{m}]\nonumber\\
&=&(\partial a')_{0}a\otimes t^{m+n}.
\end{eqnarray}
For $b\in B,\; a'\in A,\; m,n\in \Z$,
by (\ref{edef-bracket-2})--(\ref{edef-bracket-3}) we have
\begin{eqnarray}\label{ebDa}
[b\otimes t^{m}, \hat{\partial}(a'\otimes t^{n})]
&=&[b\otimes t^{m}, \partial a'\otimes t^{n}+na'\otimes t^{n-1}]\nonumber\\
&=&b_{0}\partial a'\otimes t^{m+n}+m b_{1}\partial a'\otimes t^{m+n-1}\nonumber\\
& &+nb_{0}a'\otimes t^{m+n-1}
\end{eqnarray}
and
\begin{eqnarray}\label{eDab}
[\hat{\partial}(a'\otimes t^{n}),b\otimes t^{m}]
&=&[\partial a'\otimes t^{n}+na'\otimes t^{n-1},b\otimes t^{m}]\nonumber\\
&=&(\partial a')_{0}b\otimes t^{m+n}+n(\partial a')_{1}b\otimes
t^{m+n-1}\nonumber\\
& &+na'_{0}b\otimes t^{m+n-1}.
\end{eqnarray}
If (\ref{equotient-1})--(\ref{equotient-4}) hold, then the right-hands
of (\ref{eaDa}), (\ref{eDa-a}) and (\ref{eDab}) are straightly zero
and the right hand side of (\ref{ebDa}) is contained in $\hat{\partial}L(A)$
as
\begin{eqnarray*}
& &b_{0}\partial a'\otimes t^{m+n}+m b_{1}\partial a'\otimes t^{m+n-1}
+nb_{0}a'\otimes t^{m+n-1}\nonumber\\
&=&\partial b_{0}a'\otimes t^{m+n}+(m+n)b_{0}a'\otimes t^{m+n-1}\nonumber\\
&=&\hat{\partial}(b_{0}a'\otimes t^{m+n}).
\end{eqnarray*}
This proves that the conditions
(\ref{equotient-1})--(\ref{equotient-4}) are sufficient for
$\hat{\partial}L(A)$ to be a two sided ideal.

Conversely, assume that $\hat{\partial}L(A)$ is a two-sided ideal. Let
$a,a'\in A$.  Since
$$[a\otimes t^{m}, \hat{\partial}(a'\otimes t^{n})]\in (\hat{\partial}L(A))_{(-m-n-1)}$$
for $m,n\in \Z$, from (\ref{eaDa}) we have
\begin{eqnarray*}
a_{0}\partial a'\otimes t^{m+n}=\partial a''\otimes t^{m+n+1}+(m+n+1)a''\otimes t^{m+n}
\end{eqnarray*}
for some $a''\in A$. Because $a_{0}\partial a'\in A$ and $\partial a''\in B$, we must have
\begin{eqnarray*}
\partial a''&=&0,\\
a_{0}\partial a'&=&(m+n+1)a''. \end{eqnarray*}
By taking $m=-n-1$ we get
$a_{0}\partial a'=0$, proving (\ref{equotient-1}).

By using (\ref{eDa-a}), we have
$$(\partial a')_0a\otimes t^{m+n}=\partial a''\otimes t^{m+n+1}+(m+n+1)a''\otimes
t^{m+n}$$ for some $a''\in A$. This implies that
\begin{eqnarray*}
\partial a''&=&0,\\
(\partial a')_0a&=&(m+n+1)a''.
\end{eqnarray*}
Since $m,n$ are arbitrary, we have $(\partial a')_0a=0$. This proves
(\ref{equotient-1new}).

Similarly, for $a\in A,\;b\in B$, using (\ref{ebDa}) we have
\begin{eqnarray*}
& &b_{0}\partial a'\otimes t^{m+n}+(m b_{1}\partial a' +nb_{0}a')\otimes
t^{m+n-1}\\
&=&\partial a''\otimes t^{m+n}+(m+n)a''\otimes t^{m+n-1}
\end{eqnarray*}
for some $a''\in A$. This implies that
\begin{eqnarray*}
b_{0}\partial a'&=&\partial a'',\\
mb_{1}\partial a'+nb_{0}a'&=&(m+n)a''.
\end{eqnarray*}
We immediately have $$b_{1}\partial a'=b_{0}a'=a''.$$ Thus
$b_{0}\partial a'=\partial a''=\partial b_{0}a'.$ These prove
(\ref{equotient-2})--(\ref{equotient-2more}).

Using (\ref{eDab}) we have
\begin{eqnarray*}
& &(\partial a')_{0}b\otimes t^{m+n}+n(\partial a')_{1}b\otimes t^{m+n-1}
+na'_{0}b\otimes t^{m+n-1}\nonumber\\
&=&\partial a''\otimes t^{m+n}+(m+n)a''\otimes t^{m+n-1}
\end{eqnarray*}
for some $a''\in A$. This implies that
\begin{eqnarray*}
(\partial a')_{0}b=\partial a''\;\;\;\mbox{ and }\;\;\;n(\partial a')_{1}b+na'_{0}b=(m+n)a''.
\end{eqnarray*}
Since $m,n$ are arbitrary, we must have that
$$(\partial a')_{1}b=-a'_{0}b,\ \ \mbox{and}\ \ a''=0.$$ Thus $(\partial a')_{0}b=0$.
This completes the proof.
\end{proof}

Set
\begin{eqnarray}
{\mathcal{L}}(A\oplus B)=L(A\oplus B)/\hat{\partial}L(A).
\end{eqnarray}
If the conditions (\ref{equotient-1})--(\ref{equotient-4}) hold,
in view of Proposition \ref{pto-quotient},
we have a $\Z$-graded (quotient) nonassociative algebra
${\mathcal{L}}(A\oplus B)$.

Next, we have:

\bp{pskew-symmetry} In the setting we assume that
(\ref{equotient-1})--(\ref{equotient-4}) hold.  The multiplication for
the nonassociative algebra ${\mathcal{L}}(A\oplus B)$ is
skew-symmetric if and only if
\begin{eqnarray}
a_{0}b&=&-b_{0}a,\label{esymmetry-a0b}\\
b_{1}b'&=&b'_{1}b,\label{esymmetry-b1b'}\\
b_{0}b'&=&-b'_{0}b+\partial b'_{1}b\;\;\;\mbox{  for }a\in A,\;b,b'\in B.
\label{esymmetry-b0b'}
\end{eqnarray}
\ep

\begin{proof} For $a\in A,\; b\in B$, in view of (\ref{edef-bracket-2}),
(\ref{edef-bracket-22}) and (\ref{ehatDla}),
the skew-symmetry relation
$$[a\otimes t^{m},b\otimes t^{n}]\equiv -[b\otimes t^{n},a\otimes t^{m}]
\;\;\;\mbox{ mod }\hat{\partial}L(A)$$ for $m,n\in \Z$ is
equivalent to
\begin{eqnarray}\label{e29}
(a_{0}b+b_{0}a)\otimes t^{m+n}=\partial a'\otimes
t^{m+n+1}+(m+n+1)a'\otimes t^{m+n}
\end{eqnarray}
for some $a'\in A$, and hence it amounts to that
$$\partial a'=0\;\;\mbox{ and }\; a_{0}b+b_{0}a=(m+n+1)a'.$$
By taking $m=-n-1$, we have $a_0b+b_0a=0$. This proves that
(\ref{e29}) is equivalent to (\ref{esymmetry-a0b}).

Now consider the skew symmetry relation
\begin{equation}\label{ebbskew3.41}
[b\otimes t^{m},b'\otimes t^{n}]+[b'\otimes t^{n},b\otimes t^{m}]\equiv 0
\;\;\;\mbox{ mod }\hat{\partial}L(A)
\end{equation}
for $b,b'\in B,\; m,n\in \Z$. By
(\ref{edef-bracket-3}), this is equivalent to
\begin{equation}\label{eskew-symmetry-proof}
b_{0}b'\otimes t^{m+n}+mb_{1}b'\otimes t^{m+n-1} +b'_{0}b\otimes
t^{m+n}+nb'_{1}b\otimes t^{m+n-1}\in \hat\partial L(A).
\end{equation}
Since
\begin{equation*}
(b_{0}b'+b_{0}'b)\otimes t^{m+n}+(mb_{1}b'+nb'_{1}b)\otimes t^{m+n-1}
\in L(A\oplus B)_{(-m-n)},
\end{equation*}
(\ref{eskew-symmetry-proof}) amounts to that
\begin{eqnarray}\label{eamounts}
& &(b_{0}b'+b_{0}'b)\otimes t^{m+n}+(mb_{1}b'+nb'_{1}b)\otimes
t^{m+n-1}\nonumber\\
&=&\partial a\otimes t^{m+n}+(m+n)a\otimes t^{m+n-1}
\end{eqnarray}
for some $a\in A$. The relation (\ref{eamounts}) is equivalent to that
\begin{equation}\label{eproof-38}
b_{0}b'+b_{0}'b=\partial a,\;\;\; mb_{1}b'+nb'_{1}b=(m+n)a
\end{equation}
for all $m,n\in \Z$. Clearly, (\ref{eproof-38}) amounts to that
$a=b_{1}b'=b'_{1}b$, hence $$b_{0}b'+b_{0}'b=\partial b'_{1}b.$$ This
proves that (\ref{ebbskew3.41}) is equivalent to
(\ref{esymmetry-b1b'}) and (\ref{esymmetry-b0b'}).
\end{proof}

Furthermore, we have:
\bp{plie-algebra} Under the setting, assume
that all the conditions in Propositions \ref{pto-quotient} and
\ref{pskew-symmetry} hold.  Then the quotient nonassociative
algebra ${\mathcal{L}}(A\oplus B)$ is a Lie algebra if and only if
\begin{eqnarray}
& &a_{0}u_{0}v-u_{0}a_{0}v=(a_{0}u)_{0}v,\label{elie-algebra-1}\\
& &u_{0}v_{0}w-v_{0}u_{0}w=(u_{0}v)_{0}w,\label{elie-algebra-2}\\
& &u_{0}v_{1}w-v_{1}u_{0}w=(u_{0}v)_{1}w,\label{elie-algebra-3}\\
& &u_{1}v_{0}w-v_{0}u_{1}w-(u_{0}v)_{1}w-(u_{1}v)_{0}w=0,
\label{elie-algebra-4}
\end{eqnarray}
for $a\in A,\; u,v,w\in B$.
\ep

\begin{proof} With Propositions \ref{pto-quotient}
and \ref{pskew-symmetry} we only need to show that the Jacobi
identity relations
\begin{eqnarray*}
& &[u\otimes t^{m},[v\otimes t^{n},w\otimes t^{k}]] -[v\otimes
t^{n},[u\otimes t^{m},w\otimes t^{k}]]\\
& \equiv &[[u\otimes t^{m},v\otimes t^{n}],w\otimes t^{k}]
\;\;\;\mbox{ mod }\hat{\partial}L(A)
\end{eqnarray*}
for $u,v,w\in A\cup B$ and $m,n,k\in \Z$ are equivalent to
the relations (\ref{elie-algebra-1})--(\ref{elie-algebra-4}).

Case I: If $u,v,w\in A$, by (\ref{edef-bracket-1}) all the three terms are
straightly zero in $L(A\oplus B)$.

Case II: If $u,v\in A,\; w\in B$, 
by (\ref{edef-bracket-1})--(\ref{edef-bracket-22}) all the three
terms are also straightly zero in $L(A\oplus B)$. For example,
since $v_{0}w\in A$, by definition we have
$$[u\otimes t^{m},[v\otimes t^{n},w\otimes t^{k}]]
=[u\otimes t^{m},v_{0}w\otimes t^{n+k}]=0.$$

Case III: Assume that $u\in A,\; v,w\in B$. We have
\begin{eqnarray*}
[u\otimes t^{m},[v\otimes t^{n},w\otimes t^{k}]]
&=&[u\otimes t^{m},v_{0}w\otimes t^{n+k}+nv_{1}w\otimes t^{m+k-1}]\\
&=&u_{0}v_{0}w\otimes t^{m+n+k},\\
{[v\otimes t^{n},[u\otimes t^{m},w\otimes t^{k}]]}
&=&[v\otimes t^{n},u_{0}w\otimes t^{m+k}]\\
&=&v_{0}u_{0}w\otimes t^{m+n+k},\\
{[[u\otimes t^{m},v\otimes t^{n}],w\otimes
t^{k}]}&=&[u_{0}v\otimes
t^{m+n},w\otimes t^{k}]\\
&=&(u_{0}v)_{0}w\otimes t^{m+n+k}.
\end{eqnarray*}
In this case the Jacobi identity is equivalent to
\begin{eqnarray*}
& &(u_{0}v_{0}w-v_{0}u_{0}w-(u_{0}v)_{0}w)\otimes t^{m+n+k}\\
&=&\partial a\otimes t^{m+n+k+1}+(m+n+k+1)a\otimes t^{m+n+k}
\end{eqnarray*}
for some $a\in A$. This amounts to that
$$u_{0}v_{0}w-v_{0}u_{0}w-(u_{0}v)_{0}w=(m+n+k+1)a$$ for $m,n,k\in
\Z$ and that $\partial a=0.$ Since $m,n,k$ are arbitrary, we have
$$u_{0}v_{0}w-v_{0}u_{0}w-(u_{0}v)_{0}w=0=a.$$ 
Hence, the Jacobi identity for this case is equivalent to
(\ref{elie-algebra-1}).

Case IV: Assume that $u,v,w\in B$. We have
\begin{eqnarray*}
[u\otimes t^{m},[v\otimes t^{n},w\otimes t^{k}]]
&=&[u\otimes t^{m},v_{0}w\otimes t^{n+k}+nv_{1}w\otimes t^{n+k-1}]\\
&=&u_{0}v_{0}w\otimes t^{m+n+k}+mu_{1}v_{0}w\otimes t^{m+n+k-1}\\
& & +nu_{0}v_{1}w\otimes t^{m+n+k-1},\\
{[v\otimes t^{n},[u\otimes t^{m},w\otimes t^{k}]]}
&=&v_{0}u_{0}w\otimes t^{m+n+k}+nv_{1}u_{0}w\otimes t^{m+n+k-1}\\
& &+mv_{0}u_{1}w\otimes t^{m+n+k-1},\\
{[[u\otimes t^{m},v\otimes t^{n}],w\otimes
t^{k}]}&=&(u_{0}v)_{0}w\otimes
t^{m+n+k}+(m+n)(u_{0}v)_{1}w\otimes t^{m+n+k-1}\\
& &+m(u_{1}v)_{0}w\otimes t^{m+n+k-1}.
\end{eqnarray*}
In this case the Jacobi identity amounts to that
$$u_{0}v_{0}w-v_{0}u_{0}w-(u_{0}v)_{0}w=\partial a',$$
and
\begin{eqnarray*}
& &m(u_{1}v_{0}w-v_{0}u_{1}w-(u_{0}v)_{1}w-(u_{1}v)_{0}w)
+n(u_{0}v_{1}w-v_{1}u_{0}w-(u_{0}v)_{1}w)\\
& &=(m+n+k)a' \end{eqnarray*} for some $a'\in A$. The later
relations for all $m,n,k$ amount to that $a'=0$ and
$$u_{1}v_{0}w-v_{0}u_{1}w-(u_{0}v)_{1}w-(u_{1}v)_{0}w
=u_{0}v_{1}w-v_{1}u_{0}w-(u_{0}v)_{1}w=a'=0,$$ so that the Jacobi
identity amounts to the relations
(\ref{elie-algebra-2})--(\ref{elie-algebra-4}).
\end{proof}

The following is our main result in this section:

\bt{tfewer-axioms-lie} 
Let $A\oplus B$ be a graded vector space
equipped with bilinear maps (\ref{emap-D})--(\ref{emap-b1b}). All the
relations (\ref{esymmetry-a0b})--(\ref{equotient-4}) (in
Proposition \ref{pto-quotient}),
(\ref{esymmetry-a0b})--(\ref{esymmetry-b0b'}) (in Proposition
\ref{pskew-symmetry}) and
(\ref{elie-algebra-1})--(\ref{elie-algebra-3}) (in Proposition
\ref{plie-algebra}) exactly amount to that $A\oplus B$ is a
$1$-truncated conformal algebra. 
Furthermore, if $C=A\oplus B$ is
a $1$-truncated conformal algebra, then the (quotient)
nonassociative algebra $\Lie (A\oplus B)=L(A\oplus
B)/\hat{\partial}L(A)$ equipped with the grading defined in
(\ref{eZ-grading}) is a $\Z$-graded Lie algebra. 
\et

\begin{proof} Note that the relations (\ref{esymmetry-a0b})--(\ref{esymmetry-b0b'}) 
(in Proposition \ref{pskew-symmetry}) are exactly the commutativity relations
in the definition of a $1$-truncated conformal algebra.
Noticing $(\partial a)_{1}A=0$ because $(\partial a)_{1}$ is of degree $-1$,
we see that 
the relations (\ref{equotient-1new})-(\ref{equotient-4})
and (\ref{equotient-2}) (in Proposition \ref{pto-quotient}) 
are exactly the derivation relations in the definition of 
a $1$-truncated conformal algebra.
The other relations
(\ref{equotient-1})-(\ref{equotient-2more}) 
in Proposition \ref{pto-quotient} follow from 
the relations (\ref{equotient-1new})-(\ref{equotient-4})
(in Proposition \ref{pto-quotient}) and
(\ref{esymmetry-a0b})-(\ref{esymmetry-b0b'}) (in Proposition
\ref{pskew-symmetry}), as
for $a,a'\in A, u\in B$, using the relations
(\ref{equotient-1new})-(\ref{equotient-4}) and the relations
(\ref{esymmetry-a0b})-(\ref{esymmetry-b0b'}), we have
\begin{eqnarray*}
a_0\partial a'&=&-(\partial a')_0a=0,\\
u_0\partial a&=&-(\partial a)_0u+\partial((\partial
a)_1u)=\partial((\partial a)_1u)=-\partial (a_0u)=\partial (u_0a),\\
u_1\partial a&=&(\partial a)_1u=-a_0u=u_0a.
\end{eqnarray*}

Now we prove that in the presence of the derivation relations and 
commutativity relations, the associativity relations are equivalent to
the relations (\ref{elie-algebra-1})--(\ref{elie-algebra-3}) in Proposition
\ref{plie-algebra}.
First, note that (\ref{elie-algebra-4}) follows from others, as
using (\ref{elie-algebra-3}), (\ref{esymmetry-b0b'}) and
(\ref{equotient-4}) in order, we have 
\begin{eqnarray*}
u_1v_0w-v_0u_1w&=&-(v_0u)_1w\\
&=&(u_0v)_1w-(\partial(u_1v))_1w\\
&=&(u_0v)_1w+(u_1v)_0w.
\end{eqnarray*}
By Lemma \ref{lsimplify}, the associativity relations are equivalent to
the relations (\ref{elie-algebra-1})--(\ref{elie-algebra-3}) in Proposition
\ref{plie-algebra}. This proves the first assertion and the rest follows from
Proposition \ref{plie-algebra}.
\end{proof}

For the rest of this section {\em we assume that
$A\oplus B$ is a $1$-truncated conformal algebra}.
We denote the Lie algebra by $\Lie$. 
Let $\rho$ be the natural linear map:
\begin{eqnarray}
\rho: L(A\oplus B)\rightarrow {\mathcal{L}};\; u\otimes
t^{n}\mapsto u\otimes t^{n}+\hat{\partial}L(A).
\end{eqnarray}
For $u\in A\oplus B,\; n\in \Z,$ set
$$u(n)=\rho(u\otimes t^{n})=u\otimes t^{n}+\hat{\partial}L(A)\in
{\mathcal{L}}$$ 
and then form the generating function
 $$u(x)=\sum_{n\in \Z}u(n)x^{-n-1}\in \Lie [[x,x^{-1}]].$$
For an ${\mathcal{L}}$-module $W$, we also use $u_{W}(n)$ or sometimes
just $u(n)$ for the corresponding operator on $W$ and we write
$u_{W}(x)=\sum_{n}u(n)x^{-n-1}\in (\End W)[[x,x^{-1}]]$.  Writing the
defining commutator relations
(\ref{edef-bracket-1})--(\ref{edef-bracket-3}) in terms of generating
functions we have
\begin{eqnarray*}
[a(x_{1}),\alpha(x_{2})]&=&0,\\
{[a(x_{1}),b(x_{2})]}
&=&x_{2}^{-1}\delta\left(\frac{x_{1}}{x_{2}}\right)(a_{0}b)(x_{2}),\\
{[b(x_{1}),\beta(x_{2})]}
&=&x_{2}^{-1}\delta\left(\frac{x_{1}}{x_{2}}\right)(b_{0}\beta)(x_{2})
+(b_{1}\beta)(x_{2}){\partial\over\partial x_{2}}
x_{2}^{-1}\delta\left(\frac{x_{1}}{x_{2}}\right)
\end{eqnarray*}
for $a,\alpha\in A,\; b,\beta\in B$.

\br{rdlm-vertexla} {\em The Lie algebra $\Lie$, which was constructed
in Theorem \ref{tfewer-axioms-lie}, equipped with the linear map
$\rho$ and the linear map $\partial$ from $A$ to $B$, which is considered as
a partial linear map on $A\oplus B$, is precisely a vertex Lie algebra
with base space $A\oplus B$ in the sense of [DLM]. } \er

We have (cf. [DLM]): 

\bl{lcenter}
Let ${\mathcal{L}}$ be the Lie algebra constructed in
Theorem \ref{tfewer-axioms-lie}. Then we have
\begin{eqnarray}\label{eD-derivative-property}
(\partial a)(n)=-na(n-1)\;\;\;\mbox{ for }a\in A,\; n\in \Z.
\end{eqnarray}
If $a\in \ker \partial\subset A$, then $a(n)=0$ for $n\ne -1$ and $a(-1)$
lies in the center of the Lie algebra ${\mathcal{L}}$.
\el

\begin{proof}
For $a\in A,\; n\in \Z$, since
\begin{eqnarray*}
\hat{\partial}(a\otimes t^{n})=\partial a\otimes t^{n}+na\otimes t^{n-1},
\end{eqnarray*}
we have $(\partial a)(n)+na(n-1)=0$, proving
(\ref{eD-derivative-property}). If $a\in \ker \partial \subset A$, we have
$$na(n-1)=-(\partial a)(n)=0.$$ Thus $a(m)=0$ for $m\ne -1$. Furthermore,
by (\ref{equotient-4}) we have $$a_{0}b=-(\partial a)_{1}b=0$$ for
$b\in B$. Thus $$[a\otimes t^{m},b\otimes t^{n}]=a_{0}b\otimes
t^{m+n}=0 \;\;\;\mbox{ for }b\in B,\; m,n\in \Z.$$ By the
definition (\ref{edef-bracket-1}) we also have that $$[a\otimes
t^{m},a'\otimes t^{n}]=0\;\;\;\text{ for }a'\in A,\; m,n\in \Z.$$
These prove that $a(-1)$ lies in the center.
\end{proof}

For the $\Z$-graded Lie algebra $\Lie$ we have
$$\Lie=\coprod_{n\in \Z}\Lie_{(n)},$$
where for $n\in \Z$,
\begin{eqnarray}
\Lie_{(n)}=L(A\oplus B)_{(n)}/\hat{\partial }L(A)_{(n-1)}
=(A\otimes \C t^{-n-1}+B\otimes \C t^{-n})/\hat{\partial}(A\otimes \C t^{-n}).
\end{eqnarray}
In particular,
\begin{eqnarray}
\Lie_{(0)}&=&(A\otimes \C t^{-1}+B\otimes \C)/\hat{\partial}(A\otimes \C)
=A\otimes \C t^{-1}+ B/\partial A.
\end{eqnarray}
Recall from Corollary \ref{cLie} that $B/\partial A$ is a Lie
algebra with $A$ as a module. Denote the Lie algebra $B/\partial
A$ by $\g$:
$$\g=B/\partial A.$$ Then we have the semidirect product Lie algebra
$A\oplus \g$. Note that $\Lie_{(0)}$ is a subalgebra of $\Lie$.
{}From the commutator relations
(\ref{edef-bracket-1})--(\ref{edef-bracket-3}) we have
\begin{eqnarray}
\Lie_{(0)}\cong A\oplus \g
\end{eqnarray}
(the semidirect product Lie algebra).

Set
\begin{eqnarray}
\Lie_{(\pm)}&=&\coprod_{n\ge 1}\Lie_{(\pm n)},\\
\Lie_{(\le 0)}&=&\Lie_{-}\oplus \Lie_{(0)}=\coprod_{n\ge 0}\Lie_{(-n)}.
\end{eqnarray}

An ${\mathcal{L}}$-module $W$ is said to be {\em restricted} if
for any $w\in W$ and any $u\in A\oplus B$, $u(n)w=0$
for $n$ sufficiently large. If $W$ is an $\Lie$-module on which there
is an $\N$-grading $W=\coprod_{n\in \N}W(n)$ such that
$\Lie_{(m)}W(n)\subset W(m+n)$ for $m\in \Z,\; n\in \N$, then
$W$ must be a restricted $\Lie$-module.

We also set
\begin{eqnarray}
\Lie^{\ge 0}&=&\rho ((A\oplus B)\otimes \C[t])\subset \Lie,\\
\Lie^{< 0}&=&\rho ((A\oplus B)\otimes t^{-1}\C[t^{-1}])\subset
\Lie.
\end{eqnarray}
{}From the defining commutator relations
(\ref{edef-bracket-1})--(\ref{edef-bracket-3}), $\Lie^{\ge 0}$ and
$\Lie^{<0}$ are graded subalgebras and we have that $\Lie=\Lie^{\ge
0}\oplus \Lie^{<0}$ as a vector space. 
Furthermore we have
\begin{eqnarray}\label{elneg-decomp}
\Lie^{<0}=A(-1)\oplus B(-1)\oplus B(-2)\oplus \cdots,
\end{eqnarray}
where for $n\in \Z$,
\begin{eqnarray}
A(n)=\{ a(n)\;|\; a\in A\},\ \ \ \ 
B(n)=\{ b(n)\;|\; b\in B\}\subset \Lie.
\end{eqnarray}

Consider $\C$ as the trivial $\Lie^{\ge 0}$-module and then form the
following induced module
\begin{equation}
V_{\Lie}=U(\Lie)\otimes _{U(\Lie^{\ge 0})}\C.
\end{equation}
In view of the Poincar\'{e}-Birkhoff-Witt theorem, we
have
\begin{eqnarray}\label{eva-pbw}
V_{\Lie}=U(\Lie^{<0}),
\end{eqnarray}
as a vector space, so that we may consider $A\oplus
B$ as a subspace:
\begin{eqnarray}\label{eidentificationAB}
A\oplus B\rightarrow V_{\Lie};\;\; a+b\mapsto
a(-1){\bf 1}+b(-1){\bf 1}.
\end{eqnarray}
We assign $\deg \C=0$, making
$V_{\Lie}$ an $\N$-graded $\Lie$-module:
\begin{eqnarray}\label{egrading-valie}
V_{\Lie}=\coprod_{n\in\N}(V_{\Lie})_{(n)}.
\end{eqnarray}
It follows that  $V_{\Lie}$ is a restricted $\Lie$-module. Set
$${\bf 1}=1\in V_{\Lie}.$$
The following theorem was proved in [DLM];
the first assertion also follows from a theorem of [FKRW, MP]
(cf. [LL]):

\bt{tlie-vertex-algebra} There exists a unique vertex algebra
structure on $V_{\Lie}$ with ${\bf 1}$ as the vacuum vector and with
$Y(u,x)=u_{V}(x)$ for $u\in A\oplus B$.  The vertex algebra $V_{\Lie}$
equipped with the grading (\ref{egrading-valie}) is an $\N$-graded
vertex algebra and it is generated by $A\oplus B$.  Furthermore, any
restricted ${\mathcal{L}}$-module $W$ is naturally a $V_{\Lie}$-module
with $Y_{W}(u,x)=u_{W}(x)$ for $u\in A\oplus B$. Conversely, any
$V_{\Lie}$-module $W$ is naturally a restricted $\Lie$-module with
$u_{W}(x)=Y_{W}(u,x)$ for $u\in A\oplus B$.  \et

\br{rnotations}
{\em  Note that for $v\in A\oplus B\subset V_{\Lie}$, we have
\begin{eqnarray}
v_{m}=v(m)\;\;\;\mbox{ on }\; V_{\Lie}.
\end{eqnarray}
For $a\in A$, we have
\begin{eqnarray}
\D a=a_{-2}{\bf 1}=a(-2){\bf 1}=(\partial a)_{-1}{\bf 1}=\partial a.
\end{eqnarray}}
\er
For any $A\oplus \g$-module $U$, we have an induced $\Lie$-module
\begin{eqnarray*}
\Ind_{\Lie_{(0)}}^{\Lie}(U)=U(\Lie)\otimes _{U(\Lie_{\le 0})}U,
\end{eqnarray*}
naturally an $\N$-graded $\Lie$-module generated by $U$. It is
well known fact that the equivalence classes of 
$\N$-graded simple $\Lie$-modules one-to-one
correspond to the equivalence classes of simple $\Lie_{(0)}$-modules.

\section{Vertex algebras associated with vertex algebroids}
In this section, first we follow [GMS] to associate an $\N$-graded
vertex algebra $V_{B}$ to any vertex $A$-algebroid $B$ and we show
that $(V_{B})_{(0)}$ can be naturally identified with $A$ as
a commutative associative algebra and that $(V_{B})_{(1)}$ can be
naturally identified with $B$ as a vertex $A$-algebroid.  
Then we construct and
classify graded simple modules for the vertex algebra $V_{B}$ associated
with a vertex $A$-algebroid $B$.

Let $A$ be a commutative associative algebra with identity $e$
and let $B$ be a vertex $A$-algebroid, fixed throughout this section. 
By Proposition \ref{pact}, 
$C=A\oplus B$ is naturally a
$1$-truncated conformal algebra equipped with a bilinear map
\begin{eqnarray}
A\times B\rightarrow B;\;\; (a,b)\mapsto ab,
\end{eqnarray}
such that 
\begin{eqnarray}
eb=b\;\;\;\mbox{ for }b\in B,
\end{eqnarray}
and such that for $a,a'\in A,\; b,b'\in B$,
\begin{eqnarray}
& &a'(ab)-(a'a)b=-(a_{0}b)\partial a'-(a'_{0}b)\partial
a,\label{elast-one}\\
& &b_{0}(ab')=a(b_{0}b')-(a_{0}b)b',\label{einduced-axiom-ab2}\\
& &b_{0}(aa')=a(b_{0}a')+(b_{0}a)a',\label{eB-derivation}\\
& &a_{0}(a'b)=a'(a_{0}b),\label{einduced-axiom-aa2}\\
& &a(b_{1}b')-b_{1}(ab')=(a_{0}b)_{0}b',\label{einduced-axiom-ab4}\\
& &\partial (aa')=a\partial(a')+a'\partial(a),\label{einduced-axiom-daa}.
\end{eqnarray}
Notice that from (\ref{elast-one}) we immediately have
\begin{eqnarray}
a(a'b)=a'(ab).\label{einduced-axiom-aa1}
\end{eqnarray}

\begin{rem}
{\em With $e$ being the identity element of $A$, from (\ref{einduced-axiom-daa}), 
we have $\partial e=0$, hence by Lemma \ref{lcenter}
$e(n)=0$ for $n\ne -1$ and $e(-1)$ is a central element of $\Lie$,
of degree zero.}
\end{rem}

Associated to the $1$-truncated conformal algebra $C=A\oplus B$
we have the $\N$-graded vertex algebra 
$V_{\Lie}=\coprod_{n\in \N}(V_{\Lie})_{(n)}$.
Recall that $A\oplus B$ is a generating subspace of $V_{\Lie}$, where
\begin{eqnarray}
& &A=\{ a(-1){\bf 1}\;|\; a\in A\}\subset (V_{\Lie})_{(0)},\\
& &B=\{ b(-1){\bf 1}\;|\; b\in B\}\subset (V_{\Lie})_{(1)}. 
\end{eqnarray}

\bl{lprepare}
Set
\begin{eqnarray}
E_{0}&=&{\rm span} \{ e-{\bf 1},\; a(-1)a'-aa'\;|\; a,a'\in A\}
\subset (V_{\Lie})_{(0)},\\
E_{1}&=&{\rm span} \{ a(-1)b-ab\;|\; a\in A,\; b\in B\}
\subset (V_{\Lie})_{(1)},\\
E&=&E_{0}\oplus E_{1}\\
&=&{\rm span} \{ e-{\bf 1},\; a(-1)a'-aa',\;a(-1)b-ab\;|\; a,a'\in A,\; b\in
B\}\subset V_{\Lie}.\nonumber
\end{eqnarray}
Then
\begin{eqnarray}\label{evne}
v(n)E\subset E\;\;\;\mbox{ for }v\in C=A\oplus B,\; n\ge 0.
\end{eqnarray}
Furthermore, we have
\begin{eqnarray}
& &\D E_{0}\subset E_{1},\\
& &B(-1)E_{0}\subset A(-1)E_{1}+E_{1}.\label{ebe0-ae-e}
\end{eqnarray}
\el

\begin{proof} 
For $a\in A,\; n\ge 0$, since $\wt a(n)=-n-1<0$, we have
\begin{eqnarray*}
a(n)E_{0}=0.
\end{eqnarray*}
Similarly, for $b\in B,\; n\ge 1$, we have
\begin{eqnarray*}
b(n)E_{0}=0.
\end{eqnarray*}
For $b\in B$, due to the derivation property (\ref{eB-derivation}) we have
$b(0)e=0$. Then
$$b(0)(e-{\bf 1})=b_{0}e=0.$$
For $a,a'\in A,\; b\in B$, using (\ref{eB-derivation}) 
and the commutator formula we have
\begin{eqnarray*}
& &b(0)(a(-1)a'-aa')\\
&=&b(0)a(-1)a'-b(0)(aa')\\
&=&a(-1)b(0)a'+(b(0)a)(-1)a'-(b(0)a)a'-a(b(0)a') \\
&=&\left((b(0)a)(-1)a'-(b(0)a)a'\right)+\left(a(-1)b(0)a'-a(b(0)a')\right)
\in E_{0}.
\end{eqnarray*}
Thus
$$v(n)E_{0}\subset E_{0}\;\;\;\mbox{ for }v\in A\oplus B,\; n\ge 0.$$

For $a\in
A,\; b\in B$, since $\wt a(m)=-m-1$ and $\wt b(m)=-m$ for $m\in
\Z$, we have $a(n)E_{1}=0$ for $n\ge 1$ and $b(n)E_{1}=0$ for $n\ge 2$.
For $a, a'\in A$, $b\in B$, using (\ref{einduced-axiom-aa2}) we
have
\begin{eqnarray*}
a'_{0}(a_{-1}b-ab)=a_{-1}a'_{0}b-a'_{0}(ab)=a_{-1}(a'_0b)-a(a'_{0}b)\in
E_{0}.
\end{eqnarray*}
For $a\in A$, $b,b'\in B$, 
using the commutator formula and (\ref{einduced-axiom-ab2})
we have
\begin{eqnarray*}
b'_{0}(a_{-1}b-ab)&=&a_{-1}b'_{0}b+(b'_{0}a)_{-1}b-b'_{0}(ab)\nonumber\\
&=&a_{-1}b'_{0}b-(a_{0}b')_{-1}b-b'_{0}(ab)\nonumber\\
&=&a_{-1}b'_{0}b-(a_{0}b')_{-1}b-a(b'_{0}b)+(a_{0}b')b\\
&=&\left(a_{-1}b'_{0}b-a(b'_{0}b)\right)
+\left((a_{0}b')b-(a_{0}b')_{-1}b\right)
\in E_{1}
\end{eqnarray*}
and using the commutator formula and (\ref{einduced-axiom-ab4}) we have
\begin{eqnarray*}
b'_{1}(a_{-1}b-ab)=a_{-1}(b'_{1}b)+(b'_{0}a)_{0}b-b'_{1}(ab)
=a_{-1}(b'_{1}b)-a(b'_{1}b)\in E_{0}.
\end{eqnarray*}
This proves (\ref{evne}).

Next, we show $\D E_{0}\subset E_{1}$. We have
$$\D (e-{\bf 1})=\D e(-1){\bf 1}=e(-2){\bf 1}
=(\partial e)(-1){\bf 1}=0.$$
Using the $\D$-bracket formula (\ref{ed-bracket-formula})
and (\ref{einduced-axiom-daa}) we also have
\begin{eqnarray*}
\D (a(-1)a'-aa')&=&\D(a(-1)a')-\D (aa')\\
&=&a(-1)\D a'+(\D a)(-1)a'-\partial (aa')\\
&=&a(-1)\D a'-a\D a'+a'(-1)(\D a)-a'\D a\\
&\in & E_{1}
\end{eqnarray*}
for $a,a'\in A$. Therefore $\D E_{0}\subset E_{1}$.

Finally, we prove (\ref{ebe0-ae-e}).
For $b\in B$, we have
\begin{eqnarray*}
b(-1)(e-{\bf 1})=b(-1)e-b=e(-1)b-\D e(0)b-eb=e(-1)b-eb\in E_{1},
\end{eqnarray*}
using the fact that $e(n)=0$ for $n\ne -1$.
Let $b\in B,\; a,a'\in A$. Using commutator formula, the $\D$-bracket
formula (\ref{ed-bracket-formula}), (\ref{eB-derivation}), 
(\ref{einduced-axiom-daa}), the fact that $\D E_{0}\subset E_{1}$,  
and (\ref{elast-one}) we have
\begin{eqnarray*}
& &b(-1)(a(-1)a'-aa')\\
&=&a(-1)b(-1)a'+(b_{0}a)(-2)a'-b(-1)(aa')\\
&=&a(-1)a'(-1)b-a(-1)\D a'(0)b +\D (b_{0}a)(-1)a'-(b_{0}a)(-1)\D a'-b(-1)(aa')\\
&=&a(-1)a'(-1)b-a(-1)\D a'(0)b +\D (b_{0}a)(-1)a'-(b_{0}a)(-1)\D a'\\
& &\ \ \ \ -(aa')(-1)b+\D (aa')(0)b\\
&=&a(-1)a'(-1)b-a(-1)\D a'(0)b +\D (b_{0}a)(-1)a'-(b_{0}a)(-1)\D a'\\
& &\ \ \ \ -(aa')(-1)b-\D b_{0}(aa')\\
&=&a(-1)a'(-1)b-a(-1)\D a'(0)b +\D (b_{0}a)(-1)a'-(b_{0}a)(-1)\D a'\\
& &\ \ \ \ -(aa')(-1)b-\D (a(b_{0}a')+(b_{0}a)a')\\
&=&a(-1)a'(-1)b-a(-1)\D a'(0)b +\D (b_{0}a)(-1)a'-(b_{0}a)(-1)\D a'\\
& &\ \ \ \ -(aa')(-1)b-a\D (b_{0}a')-(b_{0}a')\D a-\D ((b_{0}a)a')\\
&=&a(-1)a'(-1)b -(b_{0}a)(-1)\D a'-(aa')(-1)b-(b_{0}a')\D a\\
& &+\left(a(-1)\D b_{0}a'-a\D (b_{0}a')\right)
+\D \left((b_{0}a)(-1)a'-(b_{0}a)a'\right)\\
&\equiv& a(-1)a'(-1)b -(b_{0}a)\D a'-(aa')b-(b_{0}a')\D a\;\;\;\mbox{mod }E_{1}\\
&=&a(-1)a'(-1)b -a'(ab)\\
&=& a(-1)(a'(-1)b-a'b)+a(-1)(a'b)-a(a'b)\\
&\in& A(-1)E_{1}+E_{1}.
\end{eqnarray*}
This proves that $B(-1)E_{0}\subset A(-1)E_{1}+E_{1}$.
\end{proof}

Define
\begin{eqnarray}
I_{B}=U(\Lie)\C[\D]E,
\end{eqnarray}
an $\mathcal{L}$-submodule of $V_{\Lie}$.

For convenience, we formulate the following simple and straightforward fact 
in linear algebra:

\bl{lsimple-fact}
Let $\g$ be any Lie algebra. 
Let $\{ u_{\alpha}\}$ be any standard PBW basis of $U(\g)$ 
associated with a basis of $\g$ and let
$f$ be any map from $\{ u_{\alpha}\}$ to $U(\g)$
such that the filtration number of
$f(u_{\alpha})$ is less than that of $u_{\alpha}$.
Then $\{ u_{\alpha}-f(u_{\alpha})\}$ is also 
a basis of $U(\g)$. 
\el

\bp{pideal} The $\Lie$-submodule $I_{B}$ of $V_{\Lie}$ is a
graded ideal of $V_{\Lie}$. Furthermore, we have
\begin{eqnarray}
(V_{\Lie})_{(0)}&=&(I_{B})_{(0)}\oplus A,\\
(V_{\Lie})_{(1)}&=&(I_{B})_{(1)}\oplus B.
\end{eqnarray}
\ep

\begin{proof} Clearly, $I_{B}$ is a graded subspace of
$V_{\Lie}$. Since $A\oplus B$ generates $V_{\Lie}$ as a vertex
algebra, $I_{B}$ is a left ideal. 
Since for $v\in A\oplus B,\; n\in \Z$, 
$[\D, v_{n}]=-nv_{n-1}$ on $V_{\Lie}$,
it follows that $\D I_{B}\subset I_{B}$.
{}From [LL], $I_{B}$ is a (two-sided) ideal of $V_{\Lie}$.

Recall from (\ref{eva-pbw}) and (\ref{elneg-decomp}) that
\begin{eqnarray*}
& &V_{\Lie}=U(\Lie^{<0}),\\
& &\Lie^{<0}=A(-1)\oplus B(-1)\oplus B(-2)\oplus \cdots
\end{eqnarray*}
as graded vector spaces.
We then have
\begin{eqnarray}
& &(V_{\Lie})_{(0)}=U(A(-1)){\bf 1}=S(A(-1)){\bf 1},\\
& &(V_{\Lie})_{(1)}=U(A(-1))B(-1){\bf 1}=S(A(-1))\otimes B.\ \ \ \ 
\end{eqnarray}

Noticing that $[\D,\Lie^{\ge 0}]\subset \Lie^{\ge 0}$, 
using Lemma \ref{lprepare} we have 
$\Lie^{\ge 0}\C[\D]E\subset \C[\D]E$, so that
\begin{eqnarray}\label{eIB-exp}
I_{B}=U(\Lie^{<0})\C[\D]E.
\end{eqnarray}
{}From this we have
\begin{eqnarray}
& &(I_{B})_{(0)}=U(A(-1))E_{0},\\
& &(I_{B})_{(1)}=U(A(-1))\D E_{0}+U(A(-1))B(-1)E_{0}+U(A(-1))E_{1}.
\end{eqnarray}
Using Lemma \ref{lsimple-fact}
we have $(V_{\Lie})_{(0)}=A\oplus U(A(-1))E_{0}$.
Since $\D E_{0}\subset E_{1}$ and $B(-1)E_{0}\subset A(-1)E_{1}+E_{1}$ 
(by Lemma \ref{lprepare}), we have
$(I_{B})_{(1)}=U(A(-1))E_{1}$. It follows from
Lemma \ref{lsimple-fact} that $(V_{\Lie})_{(1)}=(I_{B})_{(1)}\oplus B$.
\end{proof}

Set
\begin{eqnarray*}
V_{B}=V_{\Lie}/I_{B},
\end{eqnarray*}
an $\N$-graded vertex algebra. 
We have (cf. [GMS]):

\bt{tA-B}
Let $B$ be a vertex $A$-algebroid and let $V_{B}$ be the 
associated $\N$-graded vertex algebra. 
We have $(V_{B})_{(0)}=A$ and $(V_{B})_{(1)}=B$ (under the linear map
$v\mapsto v(-1){\bf 1}$) and
$V_{B}$ as a vertex algebra is generated by $A\oplus B$.
Furthermore, for any $n\ge 1$,
\begin{eqnarray*}
& &(V_{B})_{(n)}\\
&=&{\rm span}\{b_1(-n_1)...b_k(-n_k){\bf 1}\;|\;b_i\in B,\;
n_1\geq n_2\geq ...\geq n_k\geq 1,\; n_1+...+n_k=n\}.
\end{eqnarray*}
\et

\begin{proof} It follows from Proposition \ref{pideal} that
$(V_{B})_{(0)}=A$ and $(V_{B})_{(1)}=B$. Since $A\oplus B$ generates
$V_{\Lie}$ as a vertex algebra, $A\oplus B$ generates $V_{B}$ as a
vertex algebra. Furthermore we have
$$e={\bf 1},\;\; a(-1)a'=aa',\;\; a(-1)b=ab$$
for $a,a'\in A,\; b\in B$.

Since $V_{\Lie}=U(\Lie^{<0}){\bf 1}$, with $V_{B}$ as a quotient
$\Lie$-module of $V_{\Lie}$ we have
\begin{eqnarray}
V_{B}=U(\Lie^{<0}){\bf 1}.
\end{eqnarray}
Set
\begin{eqnarray}
B(-)=B(-1)\oplus B(-2)\oplus \cdots \subset \Lie^{<0}.
\end{eqnarray}
This subspace $B(-)$ is a Lie subalgebra, as for $b,b'\in B,\; m,n\ge 1$,
\begin{eqnarray*}
& &[b(-m),b'(-n)]=(b_{0}b')(-m-n)-m (b_{1}b')(-m-n-1)\\
&=&(b_{0}b')(-m-n)+\frac{m}{m+n} (\partial b_{1}b')(-m-n)\\
&\in& B(-m-n).
\end{eqnarray*}
Set
$$K=U(B(-))B(-){\bf 1}\subset V_{B}.$$
If we prove that $K=\coprod_{n\ge 1}(V_{B})_{(n)}$, 
the last assertion will follow immediately,

Since $\Lie^{<0}=A(-1)\oplus B(-)$, we have
\begin{eqnarray}
V_{B}=U(B(-))U(A(-1)){\bf 1}=U(B(-))A=A\oplus U(B(-))B(-)A.
\end{eqnarray}
Now it suffices to prove $B(-)A\subset K$.

For $b\in B,\; a\in A$, we have
$$b(-1)a=a(-1)b-\D a_{0}b=ab-\D (a_{0}b)\in B.$$
For $n\ge 1$, we have
$$nb(-n-1)a=[\D,b(-n)]a=\D b(-n)a+b(-n)\D a=\D b(-n)a+b(-n)\partial a.$$
Since $[\D,B(-)]\subset B(-)$ and $\D {\bf 1}=0$, 
we have $\D K\subset K$. 
It then follows from induction that $b(-n)A\subset K$ for $b\in B,\; n\ge 1$.
This proves $B(-)A\subset K$, completing the proof.
\end{proof}

\begin{rem} 
{\em %Claim: $V_{B}=A\oplus U(B(-))B(-)$.
%(Make $A\oplus U(B(-))B(-)$ a $V_{B}$-module.)
A basis for the whole space $V_{B}$ was given in [GMS] 
by using a different method.}
\end{rem}

Next, we shall construct and classify $\N$-graded simple $V_{B}$-modules.
Since $V_{B}$ is a quotient vertex algebra of
$V_{\Lie}$ and since $A\oplus B$ generates $V_{B}$ as a vertex
algebra, in view of Theorem \ref{tlie-vertex-algebra}
we immediately have:

\bp{pconnection}
Let $W$ be a $V_{B}$-module. Then $W$ is naturally a restricted
module for the Lie algebra $\Lie$ with $v(n)$ acting as $v_{n}$ for
$v\in A\oplus B,\; n\in \Z$. Furthermore, the
set of $V_{B}$-submodules is precisely the set of
$\Lie$-submodules.
\ep

Recall that for a vertex $A$-algebroid $B$, we have a Lie $A$-algebroid
$B/A\partial A$. We have:

\bp{ptop}
Let $W=\coprod_{n\in \N}W(n)$ be an $\N$-graded $V_{B}$-module
with $W(0)\ne 0$.
Then $W(0)$ is an $A$-module with $aw=a_{-1}w$ for $a\in A,\; w\in
W(0)$ and $W(0)$ is a module for the Lie algebra $B/A(\partial A)$ with
$bw=b_{0}w$ for $b\in B,\; w\in W(0)$. Furthermore, $W(0)$ equipped
with these module structures is a module
for Lie $A$-algebroid $B/A\partial A$. If $W$ is graded simple, then
$W(0)$ is a simple module for Lie $A$-algebroid $B/A\partial A$.
\ep

\begin{proof} First, under the defined action,
$W(0)$ is an $A$-module because
$$ew=e_{-1}w={\bf 1}_{-1}w=w\;\;\;\mbox{ for }w\in W(0)$$
and 
\begin{eqnarray*}
(aa')w=(a_{-1}a')_{-1}w
=\sum_{i\ge 0}(a_{-1-i}a'_{-1+i}w+a'_{-2-i}a_{i}w)
=a_{-1}a'_{-1}w
=a(a'w)
\end{eqnarray*}
for $a,a'\in A,\;w\in W(0)$.
Second, $W(0)$ (under the defined action) 
is a module for the Lie algebra $B/A(\partial A)$, because
for $a,a'\in A$, $b,b'\in B$ and
$w\in W(0)$, we have
\begin{eqnarray*}
& &(a\partial a'
)w=(a\partial a')_{0}w=(a_{-1}\partial a')_{0}w=\sum_{i\ge
0}(a_{-1-i}(\partial a')_{i}w+(\partial a')_{-1-i}a_{i}w)=0,\\
& &b_0b'_0w-b'_0b_0w=(b_0b')_0w=[b,b']w.
\end{eqnarray*}
Third, $W(0)$ is a module for the Lie $A$-algebroid $B/A\partial A$,
because
\begin{eqnarray*}
& &(ab)w=(a_{-1}b)_{0}w=\sum_{i\ge 0}(a_{-1-i}b_{i}w+b_{-1-i}a_{i}w)=a_{-1}b_{0}w=a(bw),\\
& &b(aw)=b_{0}a_{-1}w=a_{-1}b_{0}w+(b_{0}a)_{-1}w=a(bw)+(b_{0}a)w
\end{eqnarray*}
for $a\in A,\; b\in B$ and $w\in W(0)$.

At last, assume that $W$ is graded simple. By Proposition
\ref{pconnection}, $W$ is a graded simple $\Lie$-module.
Let $N$ be any nonzero $B/A\partial A$-submodule of $W(0)$. 
Then $\Lie_{(\le 0)}N\subset N$.
Since $A\oplus B$ generates $V_{B}$ as a vertex algebra,
$U(\Lie)N$, an $\Lie$-submodule of $W$, is a nonzero $V_{B}$-submodule.
Then $W=U(\Lie)N$ as $W$ is a graded simple
$V_B$-module. Using Poincar\'{e}-Birkhoff-Witt theorem, we have
$$W=U(\Lie)N=U(\Lie_{(+)})N=N\oplus U(\Lie_{(+)})\Lie_{(+)}N,$$ 
which implies that $W(0)=N$. This proves that
$W(0)$ is a simple $B/A \partial A$-module.
\end{proof}

On the other hand,  for any Lie $A$-algebroid
$B/A(\partial A)$-module $U$, we are going to construct an $\N$-graded
$V_{B}$-module $W=\coprod_{n\in \N}W(n)$ such that $W(0)=U$
as a module for the Lie $A$-algebroid $B/A(\partial A)$.

The following simple general result will be useful in our study:

\bl{lextra}
Let $V$ be a vertex algebra and let $I$ be a (two-sided) ideal
generated by a subset $S$. Let $(W,Y_{W})$ be a $V$-module and let
$U$ be a generating subspace of $W$ as a $V$-module such that
\begin{eqnarray}
Y_{W}(v,x)u=0\;\;\;\mbox{ for }v\in S,\; u\in U.
\end{eqnarray}
Then $Y_{W}(v,x)=0$ for $v\in I$.
\el

\begin{proof} For any subset $T$ of $V$, set
\begin{eqnarray*}
{\rm Ann}_{W}(T)=\{ w\in W\;|\; Y_{W}(a,x)w=0\;\;\;\mbox{ for all }a\in T\}.
\end{eqnarray*}
{}From [LL], we have that ${\rm Ann}_{W}(T)$ is a $V$-submodule of $W$
and that ${\rm Ann}_{W}(T)={\rm Ann}_{W}(\<T\>)$, 
where $\<T\>$ denotes the (two-sided) ideal of
$V$ generated by $T$.
Applying this to our present situation we have
$$U\subset {\rm Ann}_{W}(S)={\rm Ann}_{W}(I),$$ 
where ${\rm Ann}_{W}(I)$ is a $V$-submodule of $W$.
Since $U$ generates $W$ as a $V$-module, we must have 
${\rm Ann}_{W}(I)=W$.
That is, $Y_{W}(v,x)=0$ for all $v\in I$.
\end{proof}

Applying Lemma \ref{lextra} to our situation we have:

\begin{prop}\label{p-induction}
 Let $(W,Y_{W})$ be a $V_{\Lie}$-module.
Assume that for any $a,a'\in A, b\in B$,
\begin{eqnarray}
Y_W(e,x)u&=& u\label{econ1},\\
Y_W(a(-1)a',x)u&=&Y_W(aa',x)u,\label{econ2}\\
Y_W(a(-1)b,x)u&=&Y_W(ab,x)u,\label{econ3}
\end{eqnarray}
for all $u\in U$, where $U$ is a generating subspace of $W$ as a $V$-module.
Then $W$ is naturally a $V_B$-module.
\end{prop}

\begin{proof} Recall that $V_{B}=V_{\Lie}/I_{B}$, where
$I_{B}$ is the (two-sided) ideal of $V_{\Lie}$ generated by the subset
$$S=\{ e-{\bf 1},\;a(-1)a-aa',\;a(-1)b-ab \;|\; a,a'\in A,\; b\in B\}.$$
The $V_{\Lie}$-module $W$ is naturally a $V_B$-module if and only if
${\rm Ann}_{W}(I_{B})=W$.
Now, it follows immediately from Lemma \ref{lextra} by taking
$V=V_{\Lie}$.
\end{proof}

Let $U$ be a module for the Lie algebra $\Lie_{(0)}$.
Recall that $\Lie_{(0)}=A\oplus (B/\partial A)$.
Then $U$ is an $\Lie_{(\leq 0)}$-module by letting $\Lie_{(<0)}$ 
acts trivially on $U$.
More precisely, $U$ is an $\Lie_{(\leq 0)}$-module under the following actions
\begin{eqnarray*}
a(n-1)\cdot u&=&\delta_{n,0}au,\\
b(n)\cdot u&=&\delta_{n,0}bu
\end{eqnarray*}
for $a\in A,\; b\in B,\; n\geq 0$.
Form the induced $\Lie$-module 
\begin{eqnarray}
M(U)=\Ind_{\Lie_{(\leq 0)}}^{\Lie}U=U(\Lie)\otimes_{U(\Lie_{(\leq 0)})}U.
\end{eqnarray}
Endow $U$ with degree 0, making $M(U)$ an $\N$-graded $\Lie$-module. 
It follows that $M(U)$ is a restricted $\Lie$-module.
By Theorem \ref{tlie-vertex-algebra}, $M(U)$ is 
naturally a $V_{\Lie}$-module. 
Furthermore, by Poincar\'{e}-Birkhoff-Witt theorem, we have 
$$M(U)=S(B\otimes t^{-1}\C[t^{-1}])\otimes U\text{ (as a vector space)}.$$

We set
\begin{eqnarray}
W(U)={\rm span}\{ v_{n}u\;|\; v\in E,\; n\in \Z,\; u\in U\}\subset
M(U),
\end{eqnarray}
recalling that 
$$E={\rm span}\{e-{\bf 1},\; a(-1)a'-aa',\; a(-1)b-ab\;|\; a,a'\in
A,\; b\in B\}\subset V_{\Lie}.$$
We then set
\begin{eqnarray}
M_B(U)=M(U)/U(\Lie)W(U).
\end{eqnarray}
Clearly, $W(U)$ is a graded subspace of $M(U)$, so that
$M_B(U)$ is an $\N$-graded module.

\begin{prop}\label{emodule-key}
Let $U$ be a module for the Lie algebra 
$\Lie_{(0)}$ $(=A\oplus B/\partial A)$. Then $M_B(U)$ is a $V_B$-module.
Furthermore, if $U$ is a module for the Lie $A$-algebroid $B/A\partial
A$, then $(M_{B}(U))(0)=U$.
\end{prop}

\begin{proof} Since $U(\Lie)W(U)$ is an $\Lie$-submodule of $M(U)$,
$U(\Lie)W(U)$ is a $V_{\Lie}$-submodule of $M(U)$.
Then the quotient module $M_B(U)$ is naturally a $V_{\Lie}$-module.
Denote by $\bar{U}$ the image of $U$ in the quotient module
$M_B(U)$. Then $\bar{U}$ generates $M_B(U)$ as a $V_{\Lie}$-module.
{}From definition we have $\bar{U}\subset {\rm Ann}_{M_B(U)}(S)$.
Since the ideal $I_{B}$ is generated by $S$, 
it follows from Proposition \ref{p-induction} that
$M_B(U)$ is naturally a $V_{B}$-module.

For the second assertion we must prove that
$(U(\Lie)W(U))(0)=0$.

First, we show $W(U)(0)=0$. 
From definition, $W(U)(0)$ is spanned by the vectors
\begin{eqnarray*}
(e-{\bf 1})_{-1}u,\ \ \ (a(-1)a')_{-1}u-(aa')_{-1}u,\ \ \ \ 
(a(-1)b)_{0}u-(ab)_{0}u
\end{eqnarray*}
for $a,a'\in A,\; b\in B,\; u\in U$.
Since
$$(a(-1)a')_{-1}u=a(-1)a'(-1)u=a(a'u)=(aa')u=(aa')_{-1}u,$$
 and
$$(a(-1)b)_{0}u=a(-1)b(0)u=a(bu)=(ab)u=(ab)_{0}u.$$
we have $(U(\Lie)W(U))(0)=0$.

Second, we prove $\Lie_{(\le 0)}W(U)\subset W(U)$.
Recall from Lemma \ref{lprepare} that
$v_{n}E\subset E$ for $v\in A\oplus B,\; n\ge 0$.
For $v\in C=A\oplus B,\; n\ge 0$ and for $c\in E,\; u\in U$, 
since
$$v(n)u\in U,\; v_{i}c\in E\;\;\;\mbox{ for }i\ge 0,$$
using the commutator formula we have
\begin{eqnarray*}
v(n)Y(c,x)u=Y(c,x)v(n)u
+\sum_{i\ge 0}\binom{n}{i} x^{n-i} Y(v_i c,x)u\in W(U)((x)).
\end{eqnarray*}
If $v\in A$,  we also have
\begin{eqnarray*}
v(-1)Y(c,x)u=Y(c,x)v(-1)u
+\sum_{i\ge 0}\binom{-1}{i} x^{-1-i} Y(v_i c,x)u\in W(U)((x)),
\end{eqnarray*}
as $v(-1)u\in U$. This proves $\Lie_{(\le 0)}W(U)\subset W(U)$.
Then 
$$U(\Lie)W(U)=U(\Lie_{(+)})U(\Lie_{(\le 0)})W(U)
=U(\Lie_{(+)})W(U)=W(U)+\Lie_{(+)}U(\Lie_{(+)})W(U).$$
Since the degree zero subspaces of $W(U)$ and
$\Lie_{(+)}U(\Lie_{(+)})W(U)$ are trivial, we have $(U(\Lie)W(U))(0)=0$.
\end{proof}

Next, we continue to construct graded simple $V_B$-modules. 
Let $U$ be a module for the Lie $A$-algebroid $B/A\partial A$.
There exists a unique maximal graded
$U(\mathcal{L})$-submodule $J(U)$ of $M(U)$ 
with the property that $J(U)\cap U=0$.

\begin{thm}\label{m-sim} 
Let $U$ be a module for the Lie $A$-algebroid $B/A\partial A$.
Then $L(U)=M(U)/J(U)$ is an $\N$-graded  $V_B$-module 
such that $L(U)(0)=U$ as a module for the Lie $A$-algebroid $B/A\partial A$.
Furthermore, if $U$ is a simple $B/A\partial A$-module,
$L(U)$ is a graded simple $V_{B}$-module.
\end{thm}

\begin{proof} From construction, $L(U)$ is a graded simple
module for Lie algebra $\Lie$.
By Proposition \ref{emodule-key},
$M_{B}(U)$ is an $\N$-graded $V_{B}$-module with
$(M_{B}(U))(0)=U$. It follows that $U(\Lie)W(U)\cap U=0$.
Thus $U(\Lie)W(U)\subset J(U)$. This proves that
$L(U)$ is a quotient module of $M_{B}(U)$, hence a
$V_{B}$-module.

{}From the definition of $J(U)$, for any
nonzero graded $V_{B}$-submodule
$N$ of $L(U)$ we have $0\ne N(0)\subset L(U)(0)=U$.
If $U$ is a simple $B/A\partial A$-module, then
for any nonzero graded $V_{B}$-submodule
$N$ of $L(U)$ we have $N(0)=U$, which implies $N=L(U)$
as $U$ generates $L(U)$. Thus $L(U)$ is graded simple.
\end{proof}

Furthermore we have:

\bl{lone-direction}
Let $W=\coprod_{n\in \N}W(n)$ be an $\N$-graded simple
$V_{B}$-module with $W(0)\ne 0$. Then
$W\simeq L(W(0))$.
\el

\begin{proof} By Proposition \ref{pconnection}, $W$ is an
$\Lie$-module, hence an $\N$-graded module.
It follows that $W(0)$ is an $\Lie_{(0)}$-module.
As $W$ is $\N$-graded simple, we have $U(\Lie)W(0)=W$.
Using the universal property of $M(W(0))$,
we have an $\Lie$-homomorphism $\psi$ from $M(W(0))$ onto $W$, extending
the identity map on $W(0)$. Since $\ker \psi \cap W(0)=0$, we have
$\Ker \psi \subset J(W(0))$. Then $\psi$ reduces to an
$\Lie$-homomorphism $\bar{\psi}$ from $L(W(0))$ onto $W$.
It follows from the simplicity of $L(W(0))$ (by Theorem \ref{m-sim}) that
$\bar{\psi}$ is actually an isomorphism.
\end{proof}

To summarize we have:

\bt{tclassification}
For any complete set $H$ of representatives of equivalence classes of simple
modules for the Lie $A$-algebroid $B/A\partial A$, 
$\{ L(U)\;|\; U\in H\}$
is a complete set $H$ of representatives of equivalence classes of simple
$\N$-graded simple $V_{B}$-modules.
\et

\begin{proof}
By Lemma \ref{lone-direction},
every $\N$-graded simple $V_{B}$-module is isomorphic to $L(U)$
for some simple module $U$ for the Lie $A$-algebroid $B/A\partial A$.
Furthermore, it is clear that for simple modules $U$ and $U'$ 
for the Lie $A$-algebroid $B/A\partial A$,
$L(U)\simeq L(U')$ as $\N$-graded $V_{B}$-modules if and only if 
$U\simeq U'$
as modules for the Lie $A$-algebroid $B/A\partial A$.
\end{proof}

\end{document}